\documentclass[gdrevised]{geradwp}

\PassOptionsToPackage{hyphens}{url}

\usepackage{fancyhdr,graphicx,amsmath,amssymb}
\usepackage{subfigure}
\usepackage{multirow}
\usepackage{mathtools}

\usepackage[ruled,vlined]{algorithm2e}
\usepackage{hyperref}
\usepackage{natbib}
\bibpunct[, ]{(}{)}{,}{a}{}{,}

\def\BIBand{and}

\usepackage{natbib}
 \bibpunct[, ]{(}{)}{,}{a}{}{,}%
 \def\BIBand{and}%

\usepackage{mathtools}
\usepackage[normalem]{ulem}
\newcommand{\vrp}{VRP-VCSD}%

\usepackage{pifont}
\newcommand{\xmark}{\ding{55}}%

\graphicspath{{Figures/}} 
\hypersetup{colorlinks,%
citecolor={blue}, 
urlcolor={blue},
linkcolor={blue},
breaklinks={true}
}

\GDtitre{Off-line approximate dynamic programming for the vehicle routing problem with a highly variable customer basis and stochastic demands}
\GDmois{Septembre}{September}
\GDannee{2021}
\GDnumero{53}
\GDauteursCourts{M. Dastpak, F. Errico, O. Jabali}
\GDauteursCopyright{Dastpak, Errico, Jabali}
\GDrevised{Juin 2022}

\begin{document}

\GDpageCouverture

\begin{GDpagetitre}

\begin{GDauthlist}
\GDauthitem{Mohsen Dastpak \ref{affil:ets}\GDrefsep\ref{affil:gerad}\GDrefsep\ref{affil:cirrelt}}
\GDauthitem{Fausto Errico \ref{affil:ets}\GDrefsep\ref{affil:gerad}\GDrefsep\ref{affil:cirrelt}}
\GDauthitem{Ola Jabali \ref{affil:polimi}}
\end{GDauthlist}

\begin{GDaffillist}
\GDaffilitem{affil:ets}{Department de génie de la construction, École de technologie supérieure, Montr\'eal (Qc), Canada, H3C 1K3}
\GDaffilitem{affil:gerad}{GERAD, Montr\'eal (Qc), Canada, H3T 1J4}
\GDaffilitem{affil:cirrelt}{CIRRELT, Montr\'eal (QC), Canada, H3C 3J7}
\GDaffilitem{affil:polimi}{Dipartimento di Elettronica, Informazione e Bioingegneria, Politecnico di Milano, Piazza Leonardo da Vinci 32, Milano 20133, Italy}
\end{GDaffillist}

\begin{GDemaillist}
\GDemailitem{mohsen.dastpak.1@ens.etsmtl.ca}
\GDemailitem{fausto.errico@cirrelt.ca}
\GDemailitem{ola.jabali@ens.etsmtl.ca}
\end{GDemaillist}

\end{GDpagetitre}


\GDabstracts

\begin{GDabstract}{Abstract}
We study a stochastic variant of the vehicle routing problem (VRP) arising in the context of domestic donor collection services. The problem we consider combines the following attributes. Customers requesting services are \emph{variable}, in the sense that the customers are stochastic, but are not restricted to a predefined  set, as they may appear anywhere in a given service area. Furthermore, demand volumes are also stochastic and are observed upon visiting the customer.
The objective is to maximize the expected served demands while meeting vehicle capacity and time restrictions. We call this problem the VRP with a highly Variable Customer basis and Stochastic Demands (\vrp{}).
For this problem, we first propose a Markov Decision Process (MDP) formulation representing the classical \emph{centralized} decision-making perspective where one decision-maker establishes the routes of all vehicles. While the resulting formulation turns out to be intractable, it provides us with the ground to develop a new MDP formulation, which we call \emph{partially decentralized}. In this formulation, the action-space is decomposed by vehicle. However, the decentralization is not complete as we enforce identical vehicle-specific policies, while optimizing the collective reward. 
We propose a number of strategies to reduce the dimension of the state and action spaces associated with the partially decentralized formulation. These  yield a considerably more tractable problem, which we solve via Reinforcement Learning. In particular, we develop a Q-learning algorithm  called DecQN, which features state-of-the-art acceleration techniques such as Replay Memory and Double Q Network .
We conduct a thorough computational analysis. Results show that DecQN considerably outperforms three benchmark policies. Moreover, when comparing with existing literature, we show that our approach can compete with specialized methods developed for the particular case of the \vrp{} where customer locations and expected demands are known in advance. Finally, we show that the value functions and policies obtained by DecQN can be easily embedded in Rollout algorithms, thus further improving the performances of DecQN.

\end{GDabstract}







\GDarticlestart

\section{Introduction}
More than two billion tonnes of municipal solid waste (MSW) is annually generated.
This figure is expected to increase by 70\% by 2050 \citep{Kaza2018}.
The importance of globally containing and diminishing waste has been receiving increasing attention from policymakers. For example, waste reduction is key in several targets of the UN Sustainable Development Goals, notably target 12.5 states that ``by 2030, substantially reduce waste generation through prevention, reduction, recycling and reuse" \citep{UN2021}. This target is to be measured by national recycling rates. Recycling has also been central to EU policies. According to the 2020 EU circular economy action plan, member states must recycle or prepare for reuse at least 60\% of their municipal waste by 2030 \citep{EuropeanUnion2008}. 
Alongside the efforts invested in promoting recycling, policymakers are also exploring means of avoiding waste altogether.

Several initiatives targeting waste prevention fall under the broad title of circular economy (CE). For a comprehensive overview of the research on waste management within CE, we refer the reader to \cite{Ranjbari2021}. Within this context, we focus our attention on applications related to product reuse and the rise of online platforms to facilitate it. We then discuss a specific transportation problem arising in such platforms. 

Reuse entails extending product life cycle by having the product (or its materials) used by people that are different from its original owners \citep{Fortuna2017}. One of the most studied  applications of product reuse relates to textile and fashion items. As noted by \cite{Shirvanimoghaddam2020}, the disposal nature of \emph{fast fashion} coupled with the throwaway culture is resulting in serious environmental, social and economic problems.  
Indeed, fast fashion firms may be legally forced to collect more preowned items for reuse and recycling \citep{ZanjiraniFarahani2021}. For a comprehensive review of the environmental impact resulting from textile reuse and recycling, we refer the reader to \cite{Sandin2018}. Although less studied in the literature, another important application of product reuse relates to furniture. According to the EPA, 12.2 million tons of furniture waste generated in 2017, with 80.2\% of it ending up in landfills \citep{USEPA2019}. \cite{Curran2010} analyzed the operations of nearly 400 furniture reuse organizations (FROs) in England and Wales. The term FRO reflects the primary type of household item dealt with, yet other types of items are now commonly collected, including electrical appliances and IT equipment. The study indicates that donors often only have to wait one or two days for the collection of their items.

A number of digital-based platforms have been developed to support household-based MSW reduction and CE principles.~\cite{Gu2019} study the environmental performance of MSW recycling associated with Aibolv, a WeChat applet. This mobile application allows individual users to place requests of disposing their MSW; then nearby operators are matched and dispatched to collect it. A wide variety of MSWs are considered, including kitchen tools and spent textiles.~\cite{Rosu2017} propose the Social Needs Marketplace, which is a platform for the efficient and effective provisioning of goods for vulnerable populations. The aim of the platform is to provide an efficient way to redistribute goods by directly matching individual donors and volunteer transporters with recipients. Considering the use case of furniture donations, the platform addresses logistical problems such as scheduling collection services using a volunteer network.

In this paper, we study a distribution planning problem arising in the context of Domestic Donor Collection Services (DDCSs).
Inspired by platforms such as \cite{BigBrothersBigSisters2022} and \cite{211GreaterMontreal2022}, which deploy DDCSs, we identify a number of unique modeling characteristics.
First, similar to e-commerce applications, DDCSs deal with a very large basin of service requests, potentially comprising all residential addresses in a city.
The number of customers to be served  each day is, however, limited and generally highly variable from one day to another.
Second, DDCSs require donors to classify the size of their products from a list of preset categories. Therefore, while donors estimate the size of their products, the actual product sizes are only revealed upon the arrival of the collection service. 
The majority of transported products are dry. Thus, the products may be considered homogeneous in terms of the vehicle capacity usage.
Finally, such platforms often deploy collection services with volunteers \citep{Curran2010}. Therefore, it is likely that not all products may be collected on a given day.

Another essential feature of the application under study is a certain repetition of the information process and problem parameters.  Donors' requests are registered daily, from the same basin of donors, in a given time period called the \emph{registration period}. At the end of this period, pickup operations are planned.
Vehicle operations begin afterward. Thus, in essence, the DDCSs planning problem is rather similar from one day to another, except that the set of donors requesting service changes. On the other hand, statistical information on donors' locations and demands are generally available in the form of historical data. 
We refer to this feature by the term \emph{variable customers}.

Our main objective is to develop a planning tool that can be adopted by DDCSs to program their collection activities.
Given the nature of the application, this tool should  capture the repetitive nature of the planning problem and should not require intensive daily computational effort.
To achieve this, we first introduce a variant of the vehicle routing problem compatible with the previously described DDCS characteristics, which we call the Vehicle Routing Problem with a highly Variable Customer basis and Stochastic Demands (\vrp{}).
We assume that the set of customers is variable and that customers may request service from anywhere within a given service area.
We model the evolution of the information process as follows. First, since the service requests are received during the registration period, which precedes the planning and the dispatching period,  we assume that customer locations and expected demands are known at the beginning of the planning phase. Second, the actual customer demand is observed upon visiting the customer. 
The \vrp{} then consists in dispatching a homogeneous fleet of
vehicles initially located at the depot to serve the demands of the realized set of customers. 
Vehicles are allowed to perform so-called \emph{preventive restocking}~\citep[see, for example][]{Louveaux2018}, thus enabling the possibility to return to the depot before the vehicle capacity is filled.
Furthermore, similar to~\cite{Goodson2013, Goodson2016}, we assume that the fleet operates during a limited period of the day. Thus, a portion of demands may remain unserved, and the objective is then to maximize the total served demand. Such an objective is coherent with the DDCS nature.

Markov Decision Process (MDP) formulations provide a flexible and rich modeling framework for stochastic and dynamic decision problems. Given the complexity of the problem at hand and the variable customer feature, MDP formulations are particularly suited for the \vrp{}.
Thus, we first propose an MDP model of the \vrp{}. This is a classical  centralized model, where the decision-maker has full information on the vehicles' and clients' statuses and establishes the routes for all vehicles. The resulting model is, however, intractable even for small problem instances. The main challenge comes from the multi-vehicle nature of the \vrp{} and the consequent explosion in the dimension of the state and action spaces.
Nonetheless, this model provides us with the ground to develop what we call a \emph{partially decentralized} MPD formulation for the \vrp{}, which significantly reduces the dimension of the state and action spaces.
The term decentralized is used in the literature to mean that the value function is computed for separated vehicle-specific action spaces \citep{OroojlooyJadid2019}. 
Additionally, the partially decentralized formulation proposed in this paper has two main features. The first is based on the observation that two vehicles will choose the same action if they are in the same situation. As a result, instead of searching for one policy for each vehicle, we rather seek a unique policy that is applicable to all vehicles. As demonstrated in Section \ref{sec:decent_mdp}, this is crucial for the success of the proposed algorithm. The second feature is that vehicles have access to the full state of the system, including the information about other vehicles (i.e., centralized). However, the global state is filtered via an \emph{observation function}, which selects  a subset of the information that is most relevant for a given vehicle, e.g., the set of closest customers.
We note that although it is suggestive to think of the partially decentralized formulation in terms of self-deciding vehicles, the resulting policy may be operated in a central fashion. Indeed, in our approach, vehicles decide individually considering the status of all vehicles and customers, while maximizing the demand collected by all vehicles.

Among the available algorithms to solve routing problems expressed in terms of MDPs, we distinguish between two main families: \emph{off-line} and \emph{on-line} algorithms~\citep{Ritzinger2016}. In general, off-line methods invest most of the available computing time before operations start, and use the available information to provide an easy-to-use operational policy that can be repeatedly applied to many problem instances. On the contrary, on-line algorithms concentrate most of the computing efforts during  operations taking advantage of the information revealed during operations.
Given the repetitive nature of our application, and the limited computational resources typical of DDCSs, we turn our attention to off-line methods.
To find the resulting operational policy, we resort to Reinforcement Learning, and in particular, we implement a Q-learning algorithm \citep{Watkins1992}, which is an approximate and model-free form of the Value Iteration Algorithm for stochastic dynamic programming \citep{Powell2011} and it is an off-line algorithm.
Q-learning enables estimating the value of each state-action pair via simulation.
In its simplest form, states and actions are discrete, and state-action values are represented as a lookup table. However, this paper proposes a continuous state representation and incorporates a two-layer artificial neural network to approximate the Q values.
The proposed algorithm, called DecQN, features two main key elements.
First, the challenges of variable-sized customer sets are addressed by adopting a fixed-size vector based on a particular heatmap-style encoding technique.
Second, given the proposed partially decentralized model and our choice to develop a single policy for all vehicles, our algorithm is only required to maintain a single set of Q values, which is iteratively updated according to the simulation data from all vehicles. 
Finally, we adopt several state-of-the-art techniques to improve the overall performance of our algorithm, including Replay Memory \citep{Mnih2015} and Double Q Network \citep{VanHasselt2016}. 

We conduct a thorough computational analysis. First, given that no instance set is available for the \vrp{} in the literature, we construct a new instance set. We then implement three benchmark policies to compare with the DecQN. Results show that the DecQN considerably outperforms the considered benchmark policies. To provide a comparison with existing literature, we trained our DecQN on instances of the Vehicle Routing Problem with Stochastic Demand (VRPSD), which is a particular case of the \vrp{} where the customer set and expected demands are known in advance. Results show that even though our algorithm is not tailored to exploit the in-advance knowledge of customer locations, it is competitive with the current state-of-the-art algorithm of \cite{Goodson2016}. Further computational analysis is then geared towards testing how the training phase of the DecQN can be guided towards generalization over different instances and problem dimensions, such as duration limits and stochastic variability. Finally, we show that the obtained policies can be easily employed in an \emph{on-line} algorithm, thus further improving the performance of the DecQN. 

To summarize, the contributions of the present work are as follows:
\begin{itemize}
    \item We model the DDCSs by introducing the \vrp{}, which is a variant of the VRP with a variable customer basis and stochastic demand. For this new problem, we provide an MDP formulation based on  the traditional centralized decision-making framework. 
    \item We propose a partially decentralized MDP-based formulation of the \vrp{}, which entails searching for a single policy that is the same for all vehicles. This formulation enables us to develop computationally efficient state and action space aggregation strategies. 
    
    \item We solve the \vrp{} via a state-of-the-art Q-learning algorithm, featuring recent accelerating strategies such as Replay Memory and Double Q Networks.
    This, in combination with the partially decentralized formulation, enables efficiently sharing simulated experiences among vehicles.
    
    \item  We perform extensive computational experiments showing  that the DecQN clearly outperforms three benchmarks for \vrp{}. Furthermore, when tested on instances of the VRPSD, DecQN is competitive with state-of-the-art methods specifically designed for that problem. Moreover, we show that the training phase of the DecQN can be generalized to tackle problem instances varying in terms of several dimensions, such as duration limits and stochastic variability.
    Finally, we show that the obtained off-line policies are suitable to be used as base policies in a Rollout Algorithm.
 \end{itemize}

The paper is organized as follows. We provide a detailed literature review in Section~\ref{sec:literature}. In Section~\ref{problemdescription}, we describe the \vrp{} and provide the centralized MDP formulation. We then present the partially decentralized formulation in Section~\ref{sec:decent_mdp}, describe the proposed algorithm in Section~\ref{method}, and present our computational results in Section~\ref{results}. Finally, we present our conclusions in Section~\ref{conclusion}.

\section{Literature Review} \label{sec:literature}
We first review optimization problems related to the \vrp{} in Section~\ref{sec:litvrpscd}. We then review MDP-based solution methods for stochastic routing problems in Section~\ref{sec:litmdp}.

\subsection{Related problems}\label{sec:litvrpscd}

The \vrp{} belongs to the family of Stochastic VRPs (SVRPs), which  have been extensively studied in the literature (see \cite{Oyola2018, Oyola2017} for a comprehensive review).
The \vrp{} accounts for two sources of uncertainty; the first concerns the set of customers requesting the daily service, and the second concerns their demands.
The literature addressing customer uncertainty can be classified in two main streams. 
In the first stream, a planner receives requests up until the beginning of the operations period, e.g., 8:00 AM.
As the operations begin, no more requests are accepted, and the observed set of customers is served (e.g., \cite{Bono2020}). We refer to  this form of uncertainty as \emph{static} customers. 
In the  second stream, a planner may continue receiving requests during the operations period and decides whether to accept or reject them \citep[see, for example][]{Chen2022, Ulmer2017}. We refer to  this form of uncertainty as \textit{dynamic} customers.

Among SVRPs with static customers, we further distinguish between two categories, according to the assumptions made about the set of potential customers. The first category assumes that the set of potential customers is fixed and known in advance (e.g., \cite{Gendreau1995}). In this setting, only a subset of customers is ``present'' each day. Given that the customer set is known in advance and does not change from one day to another, we denote this as the \emph{fixed} customer assumption. 
It is worth noting that several works on stochastic versions of the Traveling Salesman Problem share this modeling framework, such as \cite{Voccia2013}.
Conversely, the second category of works assumes that customers may request service from any location in a given service area~(e.g., \cite{Lin2020, Bono2020}).
As previously mentioned, we refer to this feature as variable customers.

Works focusing on demand uncertainty generally assume that the set of customers requesting service is known in advance, i.e., fixed. 
Within this stream, we identify two main problem settings. In the first, vehicles must fully serve all customer demands. Due to the uncertainty of the latter, vehicles may have to return to the depot to reload and continue the service. In such cases, the objective is to minimize the expected total costs, which includes returns to the depot (e.g., \cite{Louveaux2018}).
In the second stream, similarly to the \vrp{}, a deadline on the duration of the operations is considered.
\cite{Mendoza2016} assume a soft deadline and penalize late services. Other authors do not allow late services, and in particular in \cite{Erera2010}, additional vehicles might be dispatched if needed, while in \cite{Goodson2013, Goodson2016} not all customers might be served, and the objective is to maximize the total demand served within the deadline.

The \vrp{} belongs to the category of SVRPs with static and variable customers, and additionally considers a second layer of uncertainty on customer demands.
To the best of our knowledge, this problem setting has not been previously addressed in the literature. 

\subsection{Solution methods}\label{sec:litmdp}
MDP formulations are commonly used to model SVRPs. This section reviews solution methods for the SVRP with stochastic customers and/or demands. We label the corresponding literature as MDP-based solution methods.
We  focus on Approximate Dynamic Programming (ADP), Reinforcement Learning (RL), and Multi-agent Reinforcement Learning. For a more general discussion, the reader is referred to \cite{Soeffker2021}. Table \ref{tab:literature} summarizes the literature that is most relevant to our work. The first part of Table \ref{tab:literature} reports works with  a single-vehicle, whereas the second part reports works with  multi-vehicles. The multi-vehicle nature of the \vrp{} 
requires considerably more complex methods. Thus, our review focuses on how  large action and state spaces are addressed in  multi-vehicle problems.

We classify MDP-based methods according to the approach used in making routing decisions.
The first approach determines, at each \emph{decision epoch}, the next customer each vehicle should visit by optimizing an estimated value function~\citep[see, for example][]{Maxwell2010, Oda2018, Chen2019}. Conversely, at each decision epoch, the second approach maintains a temporary route of which only the first customer is retained. These routes are dynamically updated at each decision epoch.
The first approach provides a simple and effective way to tackle  routing problems with variable customers, such as~\cite{Oda2018, Chen2019}. Therefore, we adopt this approach in this paper.
The second approach is more suited to situation where customer locations are known in advance, such as~\cite{Goodson2013, Goodson2016}, or in dynamic customer settings, (e.g., \cite{Ulmer2020, Joe2020}.
We remark that solution methods may also combine these two strategies, as for example,~\cite{Kullman2020} where a ride-hailing problem is addressed.

\begin{table}[!hbtp]
\begin{tiny}
    \renewcommand{\arraystretch}{1.3}
    \setlength{\tabcolsep}{2.5pt}
    \begin{tabular}{lllccllllcl}\hline
        \multirow{2}{*}{Reference}  & \multirow{2}{*}{Problem}  & \multirow{2}{*}{Obj. Func.}  & \multicolumn{2}{c}{Uncert.} & \multirow{2}{*}{Fleet}  & \multirow{2}{*}{TL}  & \multirow{2}{*}{RD} & \multirow{2}{*}{Appr.}  & \multirow{2}{*}{Coor.}  &  \multirow{2}{*}{Solution method (comp.)} \\ \cline{4-5}
         &  &  & Cust. & Other &  &  &  &  &  &  \\ \hline
        \cite{Secomandi2001}    & SVRP & min TT     & F & D    & S & - & RB & - & - & RA (On.) \\ 
        \cite{Novoa2009}        & SVRP & min TT     & F & D    & S & - & RB & - & - & RA (On.) \\ 
        \cite{Peng2020}         & SVRP & min C      & V & -     & S & - & N & - & - & PG + AM (Off.) \\
        \cite{Brinkmann2019a}   & SIRP & min No. SF & F & D    & S & - & N & - & - & H (On.) \\ 
        \cite{Ulmer2017}        & DVRP & max No. SR & D & -     & S & L & RB & - & - & RA + VFA (On.-Off.) \\
        \cite{Nazari2018}       & DVRP & max SD     & D & -     & S & TW & N & - & - & PG + AM (Off.) \\ \hline
        \cite{Fan2006}          & SVRP & min TT     & F & D    & M & - & RB & DP & \xmark & RA (On.) \\ 
        \cite{Goodson2013}      & SVRP & max SD     & F & D    & M & L & RB & DP & \xmark & RA (On.) \\ 
        \cite{Goodson2016}      & SVRP & max SD     & F & D    & M & L & RB & DP & \xmark & RA (On.) \\ 
        \cite{Bono2020}         & SVRP & min C & V & T & M & TW & N & DC & CR & PG + AM (Off.) \\ 
        \cite{Lin2020}          & SVRP & min C & V & - & M & TW & N & - & - & PG (Off.) \\
        \cite{Li2021}           & SVRP & min TT & V & - & M$^1$ & TW & N & - & - & PG + AM (Off.) \\
        \cite{Brinkmann2019}    & SIRP & min No. SF & F & D    & M & - & N & DC & AP & H (On.) \\
        \cite{Joe2020}          & DVRP & min C      & D & -     & M & TW & RB & DC & AP & SA + VFA (On.-Off.) \\
        \cite{Chen2022}         & SDD & max No. SR & D & - & M$^1$ & - & RB$^2$ & - & - & Q-learning + H (On-Off.) \\
        \cite{Ulmer2020}        & SDD & max R & D & - & M & L, TW & RB & DP & \xmark & H + VFA (On.-Off.) \\ 
        \cite{Li2021a}          & PDVRP & min C & D & - & M & TW & RB & DP & AP & Q-learning (On.-Off.) \\
        \cite{Kullman2020}      & DARP & max P      & D & -     & M & L, TW & RB + N & DC & SR & Q-learning (Off.) \\
        \cite{Maxwell2010}      & DRP & min No. DS  & D & -     & M & TW & N & DC & CR & API (Off.) \\ 
        \cite{Oda2018}          & DRP & min No. RR + IT & D & - & M & TW & N & DC & \xmark & Q-learning (Off.) \\ 
        \cite{Chen2019}         & DRP & max R & D & - & M & TW & N & DC & SR & H + PG (On.-Off.) \\ 
        \cite{Li2019}           & DRP & max R & D & - & M & - & N & DC & SR & Q-learning (Off.) \\
        \hline Our research & SVRP & max SD & V & D & M & L & N & DC & CR & Q-learning (Off.) \\ \hline
    \end{tabular}
\end{tiny}
    \caption{MDP-based solution methods proposed for routing problems}
    \label{tab:literature}
    \begin{tiny}
\noindent $^1$: Heterogeneous vehicles\\
\noindent $^2$: The decision in the MDP formulation is to either reject or choose a vehicle/drone to assign that request. Then, a heuristic decides how to accommodate that request in the route.
    \begin{itemize}
        \setlength\itemsep{0.2em}
        \item Problem: Stochastic VRP (SVRP), VRP with Dynamic Requests (DVRP), Dispatching and Re-positioning Problem (DRP), Dial-A-Ride Problem (DARP), Stochastic Inventory Routing Problem (SIRP), Same-Day Delivery (SDD), Pick-up and Delivery VRP (PDVRP)
        \item Obj. Func.: The Objective Function can be min of Travel Time (min TT), Number of Delayed Services (min No. DS), Costs (min C), Number of Rejected Requests (min No. RR), Number of Service Failures (min No. SF), and Idle Time (min IT), or max of Served Demand (max SD), Profit (max P), Revenue (max R), and Number of Served Requests (max No. SR) 
        \item Uncert.:
        \begin{itemize}
            \item Cust.: Customers can be Static and Fixed (F), Static and Variable (V), or Dynamic (D)
            \item Other: Other sources of uncertainties can be Demand (D) or Travel Time (T)
        \end{itemize}
        \item Fleet: Single vehicle (S), Multiple vehicles (M), and Unlimited vehicles (U)
        \item TL: Time limits including, a duration limit for the operation (L), or Time Windows for customers (TW)
        \item RD: The routing decision is Where to go next (N) or Route-based (RB)
        \item Appr.: The approach to cope with complexities of the MDP: Decompose the problem into multiple single-vehicle problem (DP), Decentralize the action space (DC)
        \item Coord.: How decentralized vehicles are coordinated: with Collective Rewards (CR), Shaped Rewards (SR), Assignment Problem (AP), or No coordination (\xmark)
        \item Solution method (comp.): Rollout Algorithm (RA), Approximate Policy Iteration (API), Simulated Annealing (SA), Value Function Approximation (VFA), Heuristics (H), Policy Gradient (PG), Attention Mechanism (AM). comp.: indicates how/when the majority of computation is done, On-line (On) or Off-line (Off).
    \end{itemize}
\end{tiny}
\end{table}

MDP formulations are generally solved by Dynamic Programming algorithms. These, however, suffer from the so-called curse of dimensionality. In multi-vehicle SVRPs, typically featuring very large combinatorial state and decision spaces, this curse is further aggravated. Several approaches have been proposed to cope with these difficulties. 
\cite{Fan2006} and~\cite{Goodson2013} decompose the problem into multiple single-vehicle problems.~\cite{Fan2006} use a myopic heuristic technique to partition customers into $m$ subsets, one for each available vehicle. Once the subsets are established, they are kept fixed during the execution of the algorithm.~\cite{Goodson2013} refine this technique by dynamically re-partitioning customers at each decision epoch. While both approaches drastically reduce the action and state spaces, collaboration opportunities among vehicles are not explicitly enforced.

As previously mentioned, an alternative approach to reduce the curse of dimensionality is to decentralize the decision space. We recall that the term decentralize entails that the value function is computed for separated vehicle-specific action spaces rather than for the joint action space \citep{OroojlooyJadid2019}.
This approach is based on two main hypotheses: 1) once a vehicle is assigned to a destination, it is never diverted to another destination, even if this is favorable in light of new information, and 2) at any decision epoch, it is unlikely that more than one vehicle needs to be assigned to a new destination.
The simplest way to handle decentralization is to optimize the individual reward of each vehicle, regardless of the impact on the performance of other vehicles.
However, this is equivalent to ignoring cooperation among vehicles, which may lead vehicles to compete for the reward.
For example, \cite{Oda2018} study a dynamic fleet management problem for taxi dispatching and re-positioning where assignments of taxis to customers are determined to maximize the number of requests a given taxi can individually serve.

The task of achieving some level of collaboration in a decentralized framework is very challenging and it has been extensively studied in the literature of multi-agent systems \citep{OroojlooyJadid2019}.
One possible strategy is to solve, at each decision epoch, an auxiliary master problem accessing information from the value functions of each individual vehicle. For example, in the context of a bike-sharing system, \cite{Brinkmann2019} first compute the number of expected unserved requests for each individual vehicle if assigned to each of the bike stations and then deploy a greedy algorithm to solve an assignment problem minimizing the total expected number of unserved requests. For a VRP with dynamic customers, where requests are dynamically assigned to vehicles, \cite{Joe2020} develop an approximated value function to estimate the cost of each individual vehicle route and then adopt a Simulated Annealing algorithm to search for an assignment minimizing the expected total cost.
An alternative strategy is to shape the individual rewards  to account for the influence of the vehicle's actions on the overall system performance~\citep{Li2019, Chen2019, Kullman2020}. This strategy implicitly favors vehicle coordination. For example, in an order dispatching problem,~\cite{Li2019} model the reward function as a weighted sum of the prize collected by the vehicle, while accounting for potential future rewards from other vehicles.

In our work we consider an alternative strategy. On the one hand, we enforce the sequentialization of the decision process, as  in \cite{Maxwell2010,Bono2020}, i.e., we determine the action for one vehicle at a time. Additionally, we impose that the single-vehicle policy must be the same for \emph{all} vehicles. In this context, the value function represents the reward-to-go that vehicles may \emph{collectively} earn.

Although decentralizing significantly reduces the action space, the state space remains challenging.
Various techniques have been introduced to handle this issue. For example, \cite{Ulmer2015, Joe2020, Maxwell2010} replace the actual state with a set of basis functions. Alternatively, \cite{Kullman2020, Chen2019} aggregate customers' information by a grid-like discretizing of the service region.
Other works combine grid-like discretization with the assumption that vehicles cannot fully observe the current state.
For example, in a ride-sharing application, \cite{Li2019} assume that vehicles can only access trip requests information located within a given range from their current position. 
The aggregation technique we propose in this paper is a combination of the previous ones. We aggregate customers' information using a grid-like discretization, we represent a subset of customers by basis functions, and we use a vehicle-specific \emph{observation} function to return a simplified representation of the state.
Recently, new techniques involving the adoption of Deep Neural Networks have been proposed \citep[see][for example]{Bono2020, Li2021}. Although these techniques seem promising, their implementation is computationally very expensive and usually requires complex tuning procedures.

Two main methodologies have been adopted  to solve SVRPs expressed by MDP formulations.
Approximate Policy Iteration (API) methods, such as Policy Gradient, aim at directly developing the desired policy by iteratively improving an initial policy (e.g., \cite{Bono2020}). We observe that the on-line Rollout Algorithm, such as the one developed in \cite{Goodson2013}, can be viewed as a single iteration of the Policy Iteration algorithm \citep{Bertsekas2013}. 
On the other hand, Approximate Value Iteration (AVI) methods, such as Q-learning, compute a value function, from which it is straightforward to derive the corresponding policy (e.g., \cite{Chen2022, Li2021a}). We choose to implement a Q-learning algorithm.

To conclude,  MDP-based methods are promising for SVRPs. To the best of our knowledge, no existing MDP-based method in the literature can be trivially adapted to the \vrp{}. Our work fills this gap.

\section{The centralized \vrp{} } \label{problemdescription}
In Section \ref{problemdescription1}, we formally define the \vrp{}. In Section \ref{sec:cent_form}, we introduce the centralized MDP formulation of the \vrp{}.
\subsection{Problem definition}\label{problemdescription1}

We consider a service area, where customers are located along with a depot denoted by $l_0$. A set $\mathcal{V}=\{1, \dots, m \}$ of $m$ homogeneous vehicles with capacity $Q$ is initially placed at the depot. Customers may request service from any point in the service area. The set of customers requesting service $\mathcal{C}_0$ is revealed at the beginning of the operations period (e.g., every morning), with $n=|\mathcal{C}_0|$. 
Each customer $c\in \mathcal{C}_0$ is characterized by a location $l_c$, and a probabilistic demand whose expected value $\Bar{d}_c$ is assumed to be known. Accordingly, each customer is defined by a tuple of $(l_c, \Bar{d}_c)$. The customer's actual demand $\hat{d}_c$ is only realized upon the first visit. We assume that the customer's demand is splitable in the sense that, if a vehicle is unable to fully serve a visited customer, it will serve the customer as much as possible. The rest of the demand may be later served by the same vehicle or another one. The time taken to travel between locations $i$ and $j$ is assumed to be deterministic and is denoted $\tau_{i, j}$, for $i,j \in \{l_c|c\in\mathcal{C}_0\}\cup l_0$.
Vehicles begin their trips at the depot, serve a certain subset of customers, and end their routes at the depot. Furthermore, vehicles must complete their operations before a given duration limit $L$. Vehicles may return to the depot during their trip for restocking operations after leaving a customer and before visiting the next one, even if the vehicle capacity has not been filled yet. Although our problem is to pick up items, we use the term \textit{restocking} to refer to such a policy as is typical in the literature \citep{Louveaux2018}. However, in our case, restocking entails emptying the vehicle capacity. The objective is to maximize the expected total amount of served demand.

\subsection{The centralized formulation} \label{sec:cent_form}
We now provide a centralized MDP formulation of the \vrp{}. We use the term centralized to refer to the classical approach, where routes are built by one central decision-maker.
An MDP formulation generally consists of a \textit{state} which represents the current status of the system, a control \textit{action}, and a function describing the transition of the system to a new state once an action is taken and exogenous information is revealed.
In our centralized MDP formulation, the \textit{decision epoch} is any point in time at which either a vehicle visits its destination, or when new information is revealed. This includes  i) the beginning of the horizon, when the set of customers is revealed,  and ii) during operations, when at least one vehicle arrives at a customer, or when a vehicle arrives at the depot. At each decision epoch $k$, the status of the system is  represented by a state $s_k$ as follows: 
\begin{equation}\label{eq:state}
    s_k = ([(l_c, h_c, \Bar{d}_c, \hat{d}_c)]_{c\in \mathcal{C}_0}, [(l_v, a_v, q_v)]_{v\in \mathcal{V}}, t_k),
\end{equation}
where the first two groups of components refer to the state of customers and vehicles, respectively, and the last component indicates the time.
Note that although the customer and vehicle state components are naturally indexed by the decision epoch $k$, we have simplified the notation by omitting the index $k$.
The state of each customer $c\in \mathcal{C}_0$ is represented by a tuple of $(l_c, h_c, \Bar{d}_c, \hat{d}_c)$. The customer's availability $h_c$ takes value 1 if the customer is available to be served and 0 once it has been completely served.
The customer's availability $h_c$ takes value 1 if the customer is available to be served and 0 once it is assigned to a vehicle or it has been completely served. For each vehicle $v\in \mathcal{V}$, the tuple of $(l_v, a_v, q_v)$ forms its state, where $l_v, a_v$, and $q_v$ are the destination (i.e., the location the vehicle is heading to), arrival time at destination, and the current available capacity, respectively. The time of the system at decision epoch $k$ is denoted by $t_k$.

The action set is represented by an $m$-dimensional vector of $x_k=(x_k^1, ..., x_k^m)$, where $x_k^v$ indicates the next location vehicle $v$ travels to at decision epoch $k$. Let $\Bar{\mathcal{V}}_k=\{v\in\mathcal{V}|a_v=t_k \}$ be the set of active vehicles   that arrive at their destinations and are available to take action at decision epoch $k$. At decision epoch $k$, the action $x_k$ belongs to the action space $A(s_k)$ defined as:
\begin{align}
A(s_k)= & \{ x_k\in \{\mathcal{C}_0 \cup \{l_0\}\}^m:\label{eq:act0} \\
 & x_k^v=l_v, \hspace{5pt} \forall \hspace{3pt}  v\in \mathcal{V}\setminus \Bar{\mathcal{V}}_k,\label{eq:act1} \\
 & x_k^v=l_0, \hspace{5pt} \forall \hspace{3pt}  \{v\in \Bar{\mathcal{V}}_k| q_v=0\},\label{eq:act2} \\
 & x_k^v \neq x_k^{v'},  \hspace{5pt} \forall  \hspace{3pt} \{v, v' \in \mathcal{V}| v\neq v', x_k^{v'}\neq l_0 \},\label{eq:act3}  \\
 & x_k^v\in J(s_k,v) \cup \{l_0\}, \hspace{5pt}  \forall \hspace{3pt} v \in \Bar{\mathcal{V}}_k,\label{eq:act4}  \\
 & x_k^v \neq l_0, \hspace{5pt}  \forall \hspace{3pt}  \left\{v \in \Bar{\mathcal{V}}_k| \hspace{3pt} J(s_k, v) \setminus X_k^{-v} \neq \emptyset , l_v=l_0\} \right\}\label{eq:act5} ,
\end{align}
where $X_k^{-v}=\{x_k^1, ..., x_k^{v-1}, x_k^{v+1}, ..., x_k^m\}$ denotes the set of actions at time $k$ for all vehicles except vehicle $v$, and $J(s_k, v)$ the set of reachable customers by vehicle $v$ that have not been completely served. Thus, 
\begin{equation}\label{eq:feasaction}
    J(s_k, v) = \left\{c\in \mathcal{C}_0| h_c=1, \tau_{l_v,l_c} + \tau_{l_c,l_0}\leq L - t_k \right\}.
\end{equation}

In the defined action space $A(s_k)$, Condition (\ref{eq:act1}) enforces the fact that vehicles cannot be relocated when assigned to a destination. Condition (\ref{eq:act2}) implies that vehicles with no available capacity should return to the depot to replenish. Condition (\ref{eq:act3}) ensures that no two vehicles can choose the same destination unless the destination is the depot,
while Condition (\ref{eq:act4}) ensures that vehicles travel to reachable customers in $J(s_k, v)$ or to the depot. Finally, Condition (\ref{eq:act5}) entails that waiting in the depot is only allowed when there are no available customers to visit.

The information about customers' presence and demand volumes is revealed in two phases. At the first decision epoch ($k=0$), no customer is realized yet. Thus, $\displaystyle{s_0=([], [(l_0, 0, Q)]_{v\in\mathcal{V}}, t_0)}$.
At this point, a dummy action $x_0=(l_0)_{v\in\mathcal{V}}$ transfers the system to the next decision epoch ($k=1$) where the system time and the vehicle states remain the same; however, the set of realized customer locations and expected demands populates the set of customers $\mathcal{C}_0$. Consequently, we have: $\displaystyle{s_1=([(l_c, 1, \Bar{d}_c, -1)]_{c\in \mathcal{C}_0}, [(l_0, 0, Q)]_{v\in\mathcal{V}}, t_1)}.$
For each realized customer $c\in\mathcal{C}_0$, we set $\hat{d}_c$
to $-1$ at the beginning and will update it upon its first visit. Also, $h_c$ is initially set to one for all customers.

\begin{figure}
    \centering
    \footnotesize
    \includegraphics[width=0.45\textwidth]{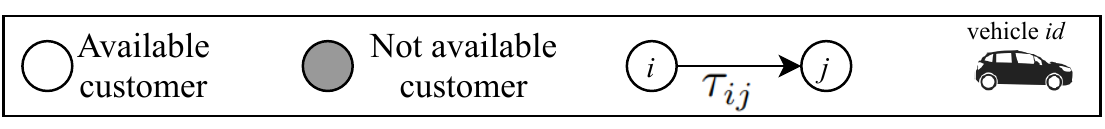}\smallskip
    
    \subfigure[No customer is realized.]{\includegraphics[width=.236\textwidth]{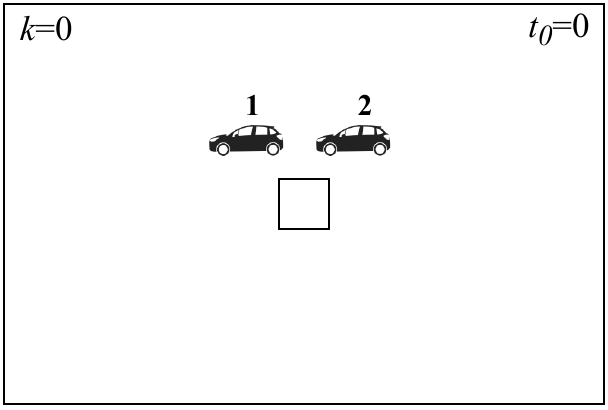}}
    ~
    \subfigure[Customers are realized.]{\includegraphics[width=.236\textwidth]{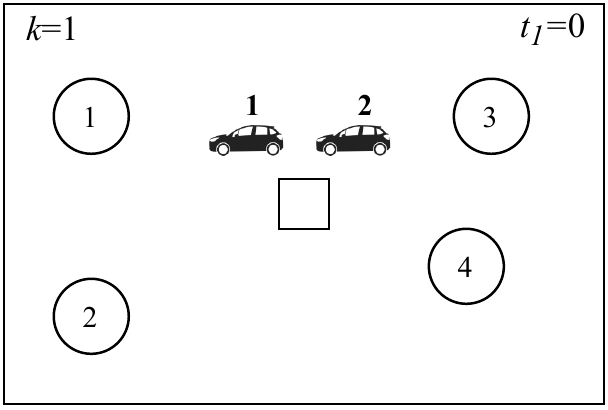}}
    ~
    \subfigure[Vehicles in $\Bar{\mathcal{V}}_1=\{1, 2\}$ choose action, $x_1=(2, 4)$. The state $s_1$ transits to the post-decision state $s_1^x$.]{\includegraphics[width=.236\textwidth]{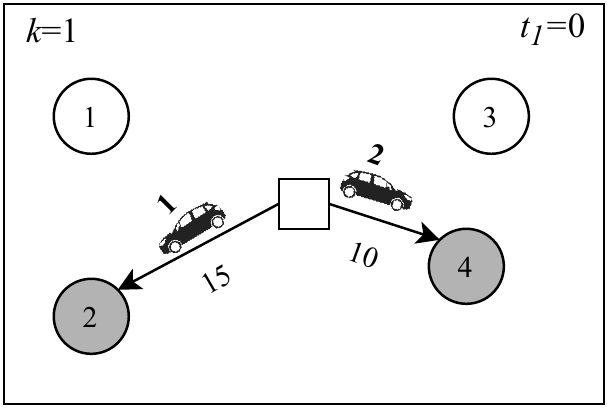}}
    ~
    \subfigure[Vehicle 2 arrives at customer 2 and observes its actual demand ($w_2^2$). The state $s^x_1$ transits to $s_2$.] {\includegraphics[width=.236\textwidth]{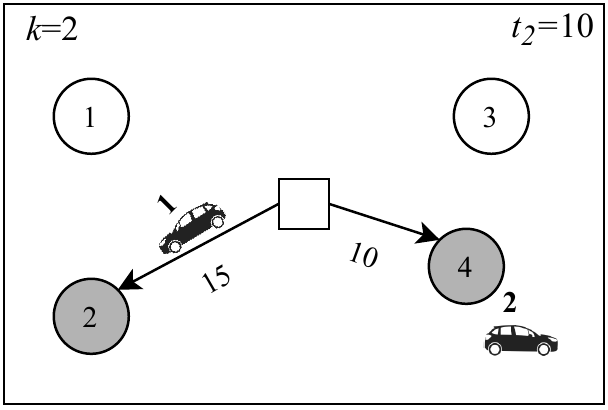}}
    \caption{The transition function in a centralized MDP formulation}
    \label{fig:trns_mdp}
\end{figure}

It is worth noting that the way the state $s_0$ is initialized highlights one of the main features of the \vrp{}. In this MDP formulation, we  develop a policy accounting for all possible customer realizations, which can tackle  realizations that have never been seen before. Accordingly, the system starts with no customers, and a random set of customers will be revealed in decision epoch 1. For $k\geq 1$, $\mathcal{C}_0$ is assumed to be fixed. Let $w$ be a sample scenario of customer demands and $w_k \subset w$ be the set of realized demand volumes of those customers whose actual demands are observed at decision epoch $k$.
By taking the action $x_k$ from the set $A(s_k)$ of possible actions in the state $s_k$, the system transits to the next state  $\displaystyle{s_{k+1} = S^M(s_k, x_k, w_{k+1})}$, where $S^M(.)$ is  a two-step transition function. In the first step, the post-decision state $s_k^x$ is the state immediately after taking action $x_k$, but before revealing the exogenous information $w_{k+1}$:
\begin{equation}
    h_{x_k^v} = 0, \hspace{5pt} \forall \hspace{3pt} \{v \in \Bar{\mathcal{V}}_k|x_k^v\neq l_0\}.
    \label{eq:pds1}
\end{equation}
\begin{equation}
    a_v = t_k + \tau_{l_v,l_{x_k^v}}, \hspace{5pt} \forall \hspace{3pt} v \in \Bar{\mathcal{V}}_k.
    \label{eq:pds2}
\end{equation}
\begin{equation}
    l_v= l_{x_k^v}, \hspace{5pt} \forall \hspace{3pt} v \in \Bar{\mathcal{V}}_k.
    \label{eq:pds3}
\end{equation}
\begin{equation}
    t_{k+1}=\min_{v\in\mathcal{V}} a_v.
    \label{eq:pds4}
\end{equation}
Equation (\ref{eq:pds1}) flags the chosen customers as unavailable. The arrival time as well as the location of active vehicles are updated using Equations (\ref{eq:pds2}) and (\ref{eq:pds3}). Finally, Equation (\ref{eq:pds4}) establishes the time of the next decision epoch, where the exogenous information $w_{k+1}$ will be revealed. All the other components are left unmodified. The second step of the transition function happens at the beginning of a new decision epoch ($k+1$) by transiting from the post-decision state $s_k^x$ to the next state $s_{k+1}$. Let $\Tilde{\mathcal{C}}_{k+1} = \{c\in\mathcal{C}_0| \exists v\in\Bar{\mathcal{V}}_{k+1} \land l_c=l_v  \}$ be the set of customers that are served
at decision epoch $k+1$.
If the actual demand of a customer in $\Tilde{\mathcal{C}}_{k+1}$ is not yet realized ($\hat{d}_c=-1$), its value will be set to the observed demand $w^c_{k+1}$.
\begin{equation}
    \hat{d}_c = w^c_{k+1}, \hspace{5pt} \forall \hspace{3pt} \{c\in \Tilde{\mathcal{C}}_{k+1}| \hat{d}_c=-1\}.
    \label{eq:trns1}
\end{equation} 
Let $\eta_v^{k+1}$ be the demand volume that vehicle $v\in\Bar{\mathcal{V}}_{k+1}$ serves at decision epoch $k+1$, defined as:
\begin{equation}
    \eta_v^{k+1} = \min \{\hat{d}_{c_v}, q_v\}, \hspace{5pt} \forall \hspace{3pt} \{v \in\Bar{\mathcal{V}}_{k+1} | l_v \neq l_0\},
    \label{eq:trns2}
\end{equation}
where $c_v$ refers to a customer $c\in\Tilde{\mathcal{C}}_{k+1}$ that is served by vehicle $v$ (i.e., $l_c=l_v$).
The unserved demands of customers in $\Tilde{\mathcal{C}}_{k+1}$ and the available capacity of vehicles in $\Bar{\mathcal{V}}_{k+1}$ are then updated as follows: 
\begin{equation}
    \hat{d}_{c_v} =\hat{d}_{c_v} - \eta_v^{k+1}, \hspace{5pt} \forall \hspace{3pt} \{v\in\Bar{\mathcal{V}}_{k+1} | l_v\neq l_0\}.
    \label{eq:trns3}
\end{equation}
\begin{equation}
    q_v = \begin{cases}
    q_v - \eta_v^{k+1} & \textrm{if  } l_v \neq l_0, \\
    Q &  \textrm{otherwise}.
    \end{cases} \hspace{5pt} \forall \hspace{3pt} v \in \Bar{\mathcal{V}}_{k+1}.
    \label{eq:trns4}
\end{equation}
Customers in $\Tilde{\mathcal{C}}_{k+1}$ whose demands are not fully served will continue being available for service. Thus,
\begin{equation}
    h_c=1, \hspace{5pt}  \forall \hspace{3pt} \{c \in \Tilde{\mathcal{C}}_{k+1}|\hat{d}_c > 0\}.
    \label{eq:trns5}
\end{equation}
These two steps are illustrated in Figure \ref{fig:trns_mdp}, which displays an example with four customers and two vehicles. 

We define the terminal decision epoch $K$ as the point in time at which all vehicles
are located at the depot, and vehicles have no other choice than to stay at the depot due to the duration limit $L$ ($A(s_K)=\{(l_0)_{v\in\mathcal{V}} \}$).
The reward function $R(s_k, x_k, w)$ is defined as the demands served by taking action $x_k$ in state $s_k$, and the exogenous information $w$ is observed:
\begin{equation} \label{eq:reward}
    R(s_k, x_k, w) = \sum_{\{v\in \Bar{\mathcal{V}}|x_k^v\neq l_0\}} \min \{q_v,  d_{x_k^v}^w\}
\end{equation}
\begin{equation}\label{eq:servdem}
    d_c^w = \begin{cases}
    w_{\bar{k}}^{c} & \textrm{if }\hat{d}_c = -1, \\
    \hat{d}_c & \textrm{otherwise}.
    \end{cases}  \hspace{5pt} \forall \hspace{3pt} c\in \mathcal{C}_0,
\end{equation}
where $w_{\bar{k}}^c$ is the demand of customer $c$ observed when the vehicle $v$ visits the customer $c$ at decision epoch $\bar{k}$.
The decision epoch $\bar{k}$ corresponds to the point in time when the vehicle $v\in\bar{\mathcal{V}}$, which has selected customer $x_k^v$ at the decision epoch $k$, arrives at its destination: 
$ \displaystyle{\bar{k}=\{k'\in [0,...,K]|t_{k'}=t_k + \tau_{l_v,l_{x_k^v}} \}}$.
Thus, the collection of the reward is generally delayed.

The value of being in state $s_k$ is defined as $V(s_k)$. This function indicates the expected cumulative reward the system can collect from state $s_k$ onward:
\begin{equation}
    V(s_k) = \max_{x_k\in A(s_k)}\mathbb{E}_{w} [R(s_k, x_k, w) + \gamma V(s_{k+1})].
    \label{eq:cent_value}
\end{equation}
In Equation \eqref{eq:cent_value}, similar to~\cite{Kullman2020}, we consider the discount factor $\gamma$ attributing greater importance to immediate rewards. 

There are some challenges making this formulation intractable for the \vrp{}. First, the dimension of the state representation is dependent on the cardinality of customers, which is stochastic. Consequently, the size of the vector $s_k$ would be variable and dependent on the number of realized customers, preventing the use of many solution methods that assume a fixed-size state. For example, no linear regression method can be implemented to approximate the value function in Equation (\ref{eq:cent_value}), when the size of the input (state vector) is variable. Furthermore, given the potentially large number of customers and the fact that the state of customers is represented by continuous variables, the state space is tremendously large.
On the other hand, managing multiple vehicles makes the size of the action space, defined in Equations (\ref{eq:act0})--(\ref{eq:act5}), excessively large. If we consider $\underline{n}$ customers available to be served at decision epoch $k$, the size of the action space is $\frac{\underline{n}!}{(\underline{n} - m)!}$. Handling such a large action set is computationally prohibitive. The following section explains how the partially decentralized formulation allows us to cope with these challenges. 

\section{The partially decentralized \vrp{}} \label{sec:decent_mdp}
The MDP formulation proposed in this section is decentralized as the action-space is decomposed by vehicle. However, the decentralization is not complete because 1) we enforce that all vehicle-specific policies are the same, 2) even though filtered by an observation function (see details in Section \ref{sec:decvrpscdstate}), we assume that vehicles access information on the global state of the system, and 3) the optimization is done over the joint policy that is obtained as the union of the vehicle-specific policies, thus the resulting value function represents the collective reward. To facilitate the decomposition of the action-space, we enforce that only one vehicle, called the \emph{active vehicle}, is allowed to make a decision at a given decision epoch.

The resulting formulation has three main  advantages. First, since the decision is restricted to one vehicle only, the size of the action space is drastically reduced to $\underline{n}$ (from potentially up to $\frac{\underline{n}}{(\underline{n}-m)!}$), where $\underline{n}$ is the number of available customers in the current decision epoch. Second, while preserving the global state of the system, the formulation allows us to filter the information that is more relevant for the active vehicle (via the observation function). 
This  substantially reduces the size of the state space, when compared to the centralized case. Third, the global policy is comprised of the union of vehicle-specific policies, which are forced to be identical. Thus, the approach has only to optimize a single policy (which is shared by all vehicles), instead of $m$ policies.

We derive the partially decentralized formulation in Section \ref{sec:decvrpscdform}.
In Section \ref{sec:decvrpscdstate}, we provide a detailed description of the vehicle observation function and its encoding. 

\subsection{The partially decentralized formulation} \label{sec:decvrpscdform}
We define two levels of state representation. The first is the \textit{global state} $s_k$, which we presented in the centralized MDP formulation. This  keeps track of the evolution of the entire system. In the second level, we introduce the observation function $O(s_k, \Bar{v})$, which takes as input the global state $s_k$ and the current active vehicle $\Bar{v}$, and returns $o_{k, \Bar{v}}$, which is the state as observed by vehicle $\bar{v}$.
The observation $o_{k, \Bar{v}}$ is comprised of four parts:
\begin{equation}
    o_{k, \Bar{v}} =[F_k|H_k|G_k|t_k],
    \label{eq:observation}
\end{equation}
where $F_k$ is the state of a subset of customers that are promising for $\Bar{v}$, $H_k$ is an aggregated overview of all customers in the service area, $G_k=[(l_v, a_v, q_v)]_{v\in\mathcal{V}}$ is the state of all vehicles, and $t_k$ is the time at decision epoch $k$. We detail the components of the observation function  in Section \ref{sec:decvrpscdstate}.

Only one vehicle at a time is allowed to make a decision. If multiple vehicles are active at decision epoch $k$, one is randomly chosen as  $\bar{v}$. We redefine the action and the action space  for $\bar{v}$ as $x_k=(x_k^{\bar{v}})$ and $A(s_k, \bar{v})$. Namely, 
\begin{align}
    A(s_k, \Bar{v})= & \{ x_k\in \{J(s_k, \Bar{v}) \cup \{l_0\}\}: \label{eq:actionspace_dec0} \\
     & x_k=l_0 \hspace{5pt} \textrm{ if } q_{\Bar{v}}=0,\label{eq:actionspace_dec1}\\
     & x_k \neq l_0 \hspace{5pt} \textrm{ if } J(s_k, \Bar{v}) \neq \emptyset, l_{\Bar{v}}=l_0 \label{eq:actionspace_dec3}\},
\end{align}
where Condition \eqref{eq:actionspace_dec0} obliges $x_k$ to be a customer from the set of feasible customers $J(s_k, \bar{v})$ or the depot, Condition (\ref{eq:actionspace_dec1}) forces the active vehicle to travel to the depot when it has no remaining capacity, and Condition (\ref{eq:actionspace_dec3}) does not allow the vehicle to stay at the depot if there are reachable customers to be served. Furthermore, the reward of taking action $x_k$ when the active vehicle $\bar{v}$ observes $o_{k, \bar{v}}$ is denoted by $R(o_{k, \bar{v}}, x_k, w)$ and remains equal to the reward function in the centralized formulation $R(s_k, x_k, w)$, if $\Bar{\mathcal{V}}_k=\{\Bar{v}\}$.

As previously mentioned, we use the global state $s_k$ to keep track of the system dynamic. Therefore, the state transition function $S^M(.)$ and its associated equations (Equations (\ref{eq:pds1})--(\ref{eq:trns5})) are still applicable here, given that $\Bar{\mathcal{V}}_k=\{\Bar{v}\}$.
Given the observation of the global state is $o_{k, \Bar{v}}$,  we aim at finding   a policy $\pi$ that indicates the best location to travel to for the active vehicle $\Bar{v}$ in the feasible action space $A(s_k, \Bar{v})$. The resulting policy selects actions such that the total served demand by all vehicles is maximized. The optimal policy $\pi^*$ is obtained by solving the following equations: 
\begin{equation}
    V^{\pi}(o_{k, \Bar{v}}) = \max_{x_k\in A(s_k, \Bar{v})}\mathbb{E}_{w} [R(o_{k,\bar{v}}, x_k, w) + \gamma V(o_{k+1, \bar{v}'})],
    \label{eq:value}
\end{equation}
\begin{equation}
    \pi^*(o_{k, \Bar{v}}) = \arg\max_{x_k\in A(s_k, \Bar{v})}\mathbb{E}_{w} [R(o_{k,\bar{v}}, x_k, w) + \gamma V(o_{k+1, \bar{v}'})],
    \label{eq:policy}
\end{equation}
where $\Bar{v}'$ is the active vehicle in decision epoch $k+1$ and $V(o_{k, \Bar{v}})$ is the expected reward-to-go starting from observation $o_{k, \Bar{v}}$ onward. In the terminal global state $s_K$, there is no customer available to be served. Hence, $\displaystyle{V^{\pi}(o_{K,v})=0, ~ \forall ~ v\in \mathcal{V}}$. 

\subsection{The observation function and its representation}\label{sec:decvrpscdstate}
The primary purpose of the observation function is to suitably select the most relevant information for the active vehicle and aggregate the remaining information. To this purpose, we identify a subset of \textit{target customers}  $\mathcal{C}_{\Bar{v}} \subset \mathcal{C}_0$ which  we deem as  more promising for the active vehicle. Therefore, instead of keeping detailed information for all customers, the observation function only maintains the comprehensive information of target customers and aggregates the rest.

Considering our problem setting, we identify two desired characteristics for the target customers  $c\in \mathcal{C}_{\Bar{v}}$: requesting large demand volumes ($d_c$), and being close to the active vehicle ($\tau_{l_c,l_{\Bar{v}}}$). 
Therefore, we define a new measure to identify $\Tilde{n}$ target customers, where $\Tilde{n}$ is an input. We choose $\Tilde{n}$ customers with higher $\displaystyle{\rho_{c, \bar{v}} = \frac{min(\Tilde{d}_c, q_{\Bar{v}})}{\tau_{l_c,l_{\Bar{v}}}}}$, with $\Tilde{d}_c=\bar{d}_c$ for $\hat{d}=-1$, and $\Tilde{d}_c=\hat{d}_c $ otherwise.
This measure normalizes the customers' demands by their distance to the active vehicle, and also accounts for the available capacity. 
The effectiveness of selecting target customers according to their score $\rho_{c, \bar{v}}$ was validated through preliminary computational tests.  

In order to maintain information related to  target customers, we represent each customer using a set of features as $=(l_c,\tau_{l_c, l_{\Bar{v}}}, \tau_{l_c, l_0}, \Tilde{d}_c, min(\Tilde{d}_c, q_{\Bar{v}}), \mu_c)$. In this feature set, $l_c$ refers to the customer's location, $\tau_{l_c, l_{\Bar{v}}}$ is the travel time between the customer and the active vehicle, $\tau_{l_c, l_0}$ is the travel time between the customer and the depot, $\Tilde{d}_c$ is the customer's unserved demand (or expected demand, if not realized yet), $min(\Tilde{d}_c, q_{\Bar{v}})$ shows the portion of customer's demand that can be served by the active vehicle, and $\mu_c$ indicates whether the actual demand is realized or not. The vector of target customers is denoted as $\displaystyle{F_k=[(l_c,\tau_{l_c, l_{\Bar{v}}}, \tau_{l_c, l_0}, \Tilde{d}_c, min(\Tilde{d}_c, q_{\Bar{v}}), \mu_c)]_{c\in\mathcal{C}_{\bar{v}}}}$.

We propose an aggregated representation of the remaining customers to provide the active vehicle with an overview of environment. The challenge here is that the stochastic customer set necessitates an aggregating technique that handles a variable number of inputs while delivering a fixed-size output.
To this purpose, we propose a heatmap-style approach to encode all customers in the service area as a fixed-size vector $H_k$.
In particular, we divide the service area into a set of square-shaped \textit{partitions} $P$. Let  $\mathcal{C}_p\subset \mathcal{C}_0$ be the subset of customers that are located in $p\in P$. We characterize each $p$ by  $(\xi^c_p, \xi^d_p)$ representing the number of customers and the total expected demand in that partition, respectively. Thus, $\displaystyle{\xi^c_p = |\mathcal{C}_p|, \hspace{10pt} \xi^d_p = \sum_{c\in \mathcal{C}_p} \Tilde{d}_c}$, and $H_k=[(\xi^c_p, \xi^d_p)]_{p\in P}$. 

In addition to $H_k$ and $F_k$, the state of the vehicles $G_k$ forms the third segment of the observation representation (Equation \eqref{eq:observation}). We define $G_k=[(l_v, a_v, q_v)]_{v\in\mathcal{V}}$, which has the same representation used in the centralized state. Therefore, the main difference between $s_k$ and $o_{k,\bar{v}}$ is the customers' state, with a variable size, replaced by a fixed-size vector $[F_k, H_k]$. Compared to the state representation $s_k$, the observation representation $o_{k,\bar{v}}$ has two advantages: it  overcomes the problem of the variable-sized state  and reduces the dimension of the observation.
In particular, the customers' state in $s_k$ is of dimension $4*|\mathcal{C}_0|$, while $[F_k, H_k]$ is of dimension $(6*\Tilde{n} + 2*|P|)$. For example, assuming  an average number of customers of $50$, the dimension of $s_k$ is $200$ on average. However, assuming  $\Tilde{n}=10$ and $|P|=25$,  the dimension of $o_{k,\bar{v}}$ will be fixed at 110.

In addition to reducing and fixing the dimension of $o_{k,\bar{v}}$, using  target customers  reduces the dimension of the active vehicle's action space by limiting its actions to target customers. Therefore, in the definition of $J(s_k, \bar{v})$, which is used to construct the action space of the active vehicle, we substitute $\mathcal{C}_0$ with $\mathcal{C}_{\Bar{v}}$. As a result, the size of the action space is reduced to $|\mathcal{C}_{\Bar{v}}|+1$, where $|\mathcal{C}_{\Bar{v}}|\leq \Tilde{n}$.

\section{Solution method} \label{method}
Both the centralized and the partially decentralized formulations suffer from the so-called curse of dimensionality, which makes them difficult to be solved by conventional stochastic dynamic programming. Approximate methods have been proposed in the literature to overcome such challenges~\citep{Powell2011}.
In this paper, we consider the Q-learning algorithm.
Our preliminary experiments suggest that Q-learning is not able to solve the centralized formulation of the \vrp{}.  This can be explained by the large number of state-action pairs entailed by this formulation. Conversely, the partially decentralized formulation, having a considerably smaller state-action space, was more suitable to be addressed by Q-learning.  We now detail the proposed algorithm for the partially decentralized formulation.

Q-learning enables learning the so-called Q factors of each state-action pair,  through   simulation. We provide a general overview of the Q-learning algorithm in Appendix A.1. The Q factor of a given state-action pair is defined as follows:
\begin{equation} \label{eq:qestimate0}
    Q(s_k, x_k) = R(s_k, x_k, w) + \gamma\max_{x_{k+1}\in A(s_{k+1})} Q(s_{k+1}, x_{k+1}).
\end{equation}
By adopting the partially decentralized formulation, we replace the state $s_k$ by the observation $o_{k,\bar{v}}$, and rewrite the Equation~\eqref{eq:qestimate0} as follows:
\begin{equation}
    Q(o_{k, \Bar{v}}, x_k) = R(o_{k, \bar{v}}, x_k, w) + \gamma\max_{x_{k+1}\in A(s_{k+1}, \Bar{v}')} Q(o_{k+1, \bar{v}'}, x_{k+1}).
    \label{eq:qcomp}
\end{equation}
We denote $\hat{Q}(o_{k, \Bar{v}}, x_k)$ as the estimation of $Q(o_{k, \Bar{v}}, x_k)$, and use it to update the Q factor of a given state-action pair by using the following equation:
\begin{equation}
    Q(o_{k, \Bar{v}}, x_k) \leftarrow Q(o_{k, \Bar{v}}, x_k) + \alpha * [\hat{Q}(o_{k, \Bar{v}}, x_k) -  Q(o_{k, \Bar{v}}, x_k)],
    \label{eq:finalqvalue}
\end{equation}
where $\Bar{v}'$ is the active vehicle at decision epoch $k+1$.

We implement the Q values in the form of an artificial neural network instead of traditional look-up tables.
This method is called Q-network in the literature, and returns the approximate Q value for a given observation-action pair ($o_{k, \bar{v}},x_k$).
For computational efficiency, it is desirable to design the Q-network in such a way that only one forward call is sufficient to return the Q values for all ($|\mathcal{C}_{\bar{v}}| + 1$) actions in $A(s_{k}, \Bar{v})$. To this purpose, we consider a neural network that returns $\Tilde{n} + 1$ Q values in every call. Notice that, under some circumstances, the number of target customers can be less than $\Tilde{n}$, making the cardinality of $A(s_{k}, \Bar{v})$ variable, but restricted to $\Tilde{n}+1$. Therefore, in the vector of $\Tilde{n}+1$ Q values approximated by the Q-network, we associate the first $|\mathcal{C}_{\bar{v}}|$ values to the Q value of target customers $\mathcal{C}_{\bar{v}}$ and assign the last one to the depot. Q values not associated with target customers nor the depot will be set to 0.

We develop a neural network with two hidden layers. This network gets the observation representation ($o_{k,\bar{v}}$) as the input layer and returns approximate Q values in the output layer. We denote the size of the input and output layers by $h^{in}$ and $h^{out}$, where $h^{in}=6\Tilde{n}+2|P|+3m+1$ and $h^{out} = \Tilde{n} + 1$. The hidden layers are assumed to be fully connected dense layers linked by an activation function (a rectified linear unit (ReLU)). The size of the hidden layers is chosen based on empirical tests. In particular, we set their sizes to $\lfloor\frac{2}{3}(h^{in}-h^{out})\rfloor + h^{out}$ and $\lfloor\frac{1}{3}(h^{in}-h^{out})\rfloor + h^{out}$, respectively.
Figure 3 in Appendix A.2 illustrates the architecture of the proposed artificial neural network.
Denoting by $\theta$ the trainable weights of the proposed neural network, the Q value of an observation-action pair is represented as $Q(o_{k, \Bar{v}}, x_k, \theta)$.
It is worth noting that any equation involving $Q(o_{k, \Bar{v}}, x_k)$ will also apply to $Q(o_{k, \Bar{v}}, x_k, \theta)$. 

To train our network, we minimize the so-called loss function. We define this function based on the difference between the new estimation of the values $\hat{Q}(o_{k, \Bar{v}}, x_k, \theta)$ and the current estimation $Q(o_{k, \Bar{v}}, x_k, \theta)$, denoted as $\Delta$ (i.e., $\Delta=\hat{Q}(o_{k,\bar{v}}, x_k, \theta) - Q(o_{k,\bar{v}}, x_k, \theta)$).
Among several loss functions used in the literature, we found that the Huber loss function performs better than the commonly used functions such as Mean Squared Error (MSE) and Mean Absolute Error (MAE). The Huber loss acts as an MSE for small $\Delta$ and as  an MAE for larger $\Delta$. The Huber loss is defined as:
\begin{equation}
    Huber(\Delta)=\begin{cases}
               \frac{1}{2}\Delta^2& \textrm{for } |\Delta| \leq \delta\\
               \delta(|\Delta| - \frac{1}{2}\delta)& \textrm{otherwise},
            \end{cases}
\end{equation}
where $\delta$ is a control parameter. We obtained the best results by setting $\delta=5$ and used the Adam optimizer \citep{Kingma2014} to update the network in order to minimize the loss function. 

Algorithm 2 (in Appendix A.2) summarizes our DecQN algorithm. Similar to what was suggested by \cite{Mnih2015}, our preliminary results showed that using the target network \citep{VanHasselt2016} and the replay memory accelerates convergence. By adopting the target network (parameterized by $\bar{\theta}$ in Algorithm 2), the model learns from a stable target. The target network $\bar{\theta}$ is periodically copied from the primary network $\theta$ every $\beta_d$ trials. In the replay memory, experiences are stored in a First-In-First-Out list, denoted as $B$, with a fixed size. At every decision epoch, with a probability of $\beta_{t}$, a random batch of experiences ($\Tilde{B}$) is sampled from this list to update the weights ($\theta$) in the neural network. Several advantages have been mentioned for using the replay memory \citep{Mnih2015} such as breaking the correlation between consecutive experiences and giving more chances to those experiences that happen more. The exploration rate ($\epsilon$) and the learning rate ($\alpha$) are linearly decayed during the training phase. The decaying rate is decided based on the length of the training process.

As previously mentioned, one of the main features of the proposed approach is that we only develop a single set of Q values, resulting in a policy shared among all vehicles. This has several advantages: 
On the one hand, the information learned by one vehicle via its interaction with the environment is shared with all other vehicles, thus accelerating the convergence of the algorithm. On the other hand, the resulting method is more sample efficient.
Finally, maintaining and storing one policy is more manageable than maintaining $m$ policies. 

\section{Experimental results} \label{results}
In this section, we present our computational experiments. Since no instance set is available in the literature, we construct a new set of instances for the \vrp{}.
We design three main experiments. The first one, presented in Section \ref{sec:resultsvrpscd1}, is aimed at evaluating the overall performance of the DecQN algorithm. 
To this purpose, we implement three benchmarks, which we call Deterministic A Priori Policy (DAPP), Random Policy (RP), and  Greedy Policy (GP). Given that DAPP requires solving a very complex deterministic optimization problem, we only test it on relatively small instances, and perform experiments on larger instances for the two other benchmarks. 
The purpose of the second experiment, presented in Section \ref{sec:resultsvrpsd1}, is to compare the DecQN with the existing literature. The challenge is that no direct comparison is possible. The closest problem in the literature we are aware of, is the VRPSD for which  \cite{Goodson2016} developed the current state-of-the-art method. As mentioned earlier, the VRPSD can be seen as a particular case of the \vrp{} where both customer locations and expected demands are known in advance.
Even though our algorithm does not take advantage of the specific characteristics of the problem, it turns out that our algorithm is competitive and sometimes outperforms the results in \cite{Goodson2016}. 
The third experiment, in Section \ref{sec:resultsaggr1}, can be seen as an extension of the second one, where we investigate the possibility of modifying the training phase of the DecQN to develop a policy able to address a larger range of VRPSD instances aggregated over specific parameters, such as the duration limit and stochastic variability.
Finally, an additional experiment providing a simple demonstration of the fact that the policies provided by our algorithm can be easily embedded as base policies in rollout algorithms, is presented in Appendix D.
We note that depending on the size of the problem, training a policy according to Algorithm 2 takes 24 to 36 hours. Problems with fewer vehicles and customers or shorter duration limits, for example, demand less training time.
The values of the hyper-parameters in Algorithm 2 and the computational environment are described in Appendix B.2.

\subsection{Performance of DecQN in \vrp{}} \label{sec:resultsvrpscd1}
We describe the main features of the instance set in Section \ref{sec:resultsvrpscdinstance}, while Appendix B.1.1 provides the details of the instance generation method. The benchmarks policies DAPP, RP and GP are described in Section \ref{sec:resultsvrpscdbench}. We report  the computational results in Section \ref{sec:resultsvrpscdres}.

\subsubsection{\vrp{} Instances}\label{sec:resultsvrpscdinstance}
As custom when building computational experiments in learning-based methods, we need to provide two types of input data, one to be used in the training phase of the algorithm, the other in the testing phase. Elements of these two sets are samples of the problem at hand. Although training and testing sets are distinct, they are typically generated by sampling from the same distributions.

We define an instance $i$ as the tuple $i=(\mathcal{A}, \Psi_{n_z}, \Psi_l, \Psi_{\bar{d}}, \Psi_{\hat{d}}, m, Q, L)$, where $\mathcal{A}$ denotes the service area partitioned into a set of zones $Z$, while $\Psi_{n_z}, \Psi_l, \Psi_{\bar{d}}, \Psi_{\hat{d}}$ denote the probability distribution functions of the number of customers $n_z$ per zone $z\in Z$ (customer density), customer locations, expected demands, and actual demands, respectively.
Furthermore, $m, Q, L$ indicate the number of vehicles, their capacity, and the duration limit, respectively.
We detail the instance generation procedure in Appendix B.1.1. In particular,  we consider five different distributions for the customer density $\Psi_{n_z}$ (Very Low, Low, Moderate, High, and Very High), and three different values of the vehicle capacity $Q$ (25, 50, and 75), resulting in a set $I$ of 15 instances.
Furthermore, the distribution $\Psi_{n_z}$  determines the values of $m$ and $L$ (Appendix B.1.1).
The resulting values of the tuple $(\Psi_{n_z},\bar{n}, m, L)$ are (Very Low, 10, 2, 143.71), (Low, 15, 2, 201.38), (Moderate, 23, 3, 221.47), (High, 53, 7, 195.54), and (Very High, 83, 11, 187.29).
Given an instance $i$, a realization $\hat{i}$ (also referred to as a sample, or a scenario) is obtained by sampling from distributions $\Psi_{n_z}, \Psi_l, \Psi_{\bar{d}}, \Psi_{\hat{d}}$. The set of all realizations $\hat{i}$ of $i$ is denoted $\hat{I}(i)$. 

\subsubsection{Benchmarks} \label{sec:resultsvrpscdbench}
We now describe the benchmark policies.
The DAPP assumes the knowledge of the customer locations and expected demands.
It then constructs a set of routes, which are evaluated in the \vrp{} setting.
To construct such routes, we formulate the deterministic counterpart of the \vrp{} (i.e., with realized customers and demands) by the MILP described in Appendix B.3. We find an initial set of routes by solving this formulation for a given set of customers with demands corresponding to their expected values. We then evaluate the performance of these routes by  simulating the classical recourse policy (i.e., once the vehicle capacity is exceeded a return trip to the depot is performed and the capacity of the vehicle is restored) for 100 demand realizations.
The MILP was solved by CPLEX 20.1 with a time limit of two hours. Unlike the DecQN, which can be instantaneously executed for any set of customers (generated from a given distribution function), this benchmark requires a significant amount of time to find the fixed routes for each customer realization.

In the RP, an action is randomly picked from the action set. We compare with this policy since it is the policy that DecQN starts with when the exploration rate in the $\epsilon$-greedy method is set to one. Therefore, the improvement of our model over RP demonstrates the effect of learning.

The GP chooses the customer $c\in A(s_k, \bar{v})$ with the highest $\Tilde{d}_c$ at each decision epoch. We recall that $\Tilde{d}_c$ provides the expected or the remaining demand of customer $c$. If several customers have the same $\Tilde{d}_c$, the closer one is selected. As GP selects customers according to the potential reward in terms of served demands, it is thus aligned with the \vrp{} objective.

RP and GP are more comparable with our DecQN. However, they cannot handle preemptive restocking. Therefore, the action set is redefined as:
\begin{equation}
    A(s_k, \bar{v}) = 
    \begin{cases}
    \{l_0\} & q_{\bar{v}} = 0 \hspace{5pt} || \hspace{5pt} J(s_k, \bar{v})=\emptyset\\
    J(s_k, \bar{v}) & \textrm{otherwise}.
    \end{cases}
\end{equation}

\subsubsection{Results and Discussion}\label{sec:resultsvrpscdres}
As previously mentioned, the benchmark DAPP requires to solve a complex deterministic MILP for each customer realization. Given the two hours time limit we impose, the DAPP was only able to address instances with Very Low and Low densities. Therefore, we organize our discussion into two parts. In the first, we present instances with Very Low and Low densities, and compare DecQN with all three benchmarks. In this case, we also restrict the vehicle capacity to 25 and 50, given that for $Q=75$ almost all demands could be served. 
In the second part, we compare the performance of DecQN with GP and the RP on instances with Moderate, High, and Very High densities.

The method to obtain the training and test sets is the same for each instance $i \in I$. In particular, the training set is obtained by sampling 5 million realizations $\hat{i}$ from $\hat{I}(i)$. We then trained the DecQN on these realizations. The exploration rate decayed linearly from 1.0 to 0.1 in the first one million trials. The learning rate also decreased linearly from $10^{-3}$ to $10^{-4}$ in the first two million trials.
 As detailed below, due to computational limitations, the cardinality of the test for the first part is smaller than the test set of the second part. Note that although the training and the test sets are sampled from the same distributions, they are disjoint.

For each instance with Very Low and Low densities, the test set is obtained by generating ten customer realizations. For each customer realization, we then generate 100 demand realizations, resulting in 1000 scenarios.
Table \ref{tab:resultsvrpscd0} shows the comparison based on average performance of DecQN, DAPP, RP, and GP. Extended results are presented in Table 11 in Appendix C. 
The first two columns in Table \ref{tab:resultsvrpscd0} indicate the characteristics of the instance, while columns $n$ and $\mathbb{E}\sum d_c$ show the average number of customers and the average total expected demand. Column Opt. reports the number of customer realizations solved to optimality by the DAPP for each instance. The average total served demand by our proposed policy is reported in column DecQN. 
Column $\%\mathbb{E}\sum d_c$ shows the portion of all demands served by DecQN. Finally, the performance gap between DecQN and the three benchmarks is reported in columns \%DAPP, \%RP, and \%GP, where measure \%X, is computed as: $\displaystyle{\textrm{\%X} = \frac{\textrm{DecQN} - X}{X}\times100}$.

The DAPP obtains optimal solutions in 17 out of 20 customer realizations on instances with Very Low density. This number decreases to 4 out of 20 for instances with Low density. Therefore, the values in Table \ref{tab:resultsvrpscd0} refer to realizations for which DAPP obtained optimal routes.  DecQN outperforms DAPP by an average of 3.0\%. 
This suggests that, when considering the Very Low and Low densities, even by giving the DAPP sufficient time to solve the resulting deterministic problem, our method is still a better choice.
Furthermore, DecQN outperforms RP and GP by, an average, 41.5\% and 19.9\%. We note that such improvements increase for $\%$RP and decrease for $\%$GP with the density.
This can be explained by the fact that, as the number of customers increases, the probability that RP deviates from the optimal solution increases. A denser instance, in contrast, is more favorable for GP. Intuitively, as instances become denser, the knapsack component of the problem becomes dominant, and greedy policies are known to behave well on these kind of problems.

According to Table \ref{tab:resultsvrpscd0}, as the vehicle capacity becomes larger (from 25 to 50 in the Very Low density), the improvement over DAPP decreases (from 2.6\% to 1.2\%). To explain this finding, we remark that the vehicle capacity directly impacts the probability of route failure. When the capacity is larger, the routes are less vulnerable to failure caused by stochastic demands.
Therefore, we argue that the advantage of our method over DAPP is more pronounced on instances with smaller capacity because our method can better handle route failures. Validating this finding for instances with Low customer density is not easy because only a few could be solved optimally.

\begin{table}[!hbtp]
\centering
\footnotesize
\begin{tabular}{cccccccccc} \hline
$\Psi_{n_z}$ & $Q$ & $n$ & Opt. & $\mathbb{E}\sum d_c$ & DecQN & \%$\mathbb{E}\sum d_c$ & \%DAPP & \%RP & \%GP \\ \hline
\multirow{2}{*}{Very Low} & 25 & \multirow{2}{*}{10.4} & \multirow{2}{*}{17/20} & \multirow{2}{*}{105.0} & 68.9 & 66.3\% & 2.6\%  & 34.8\% & 25.7\% \\
 & 50 & &  & & 87.3 & 83.5\% & 1.2\%  & 44.3\% & 15.7\% \\ \hline
\multicolumn{7}{r}{\textbf{Avg}} & \textbf{1.8\%} & \textbf{40.0\%} & \textbf{19.9\%} \\ \hline
\multirow{2}{*}{Low}  & 25 & \multirow{2}{*}{15.0} & \multirow{2}{*}{4/20} & \multirow{2}{*}{147.5} & 99.6 & 69.1\% & 2.3\%  & 58.9\% & 14.3\% \\
 & 50 & &  & & 122.0 & 84.2\% & 10.5\% & 54.5\% & 16.9\% \\ \hline
\multicolumn{7}{r}{\textbf{Avg}} & \textbf{8.5\%} & \textbf{48.1\%} & \textbf{19.8\%} \\ \hline
\multicolumn{7}{r}{\textbf{Total Avg}} & \textbf{3.0\%} & \textbf{41.5\%} & \textbf{19.9\%} \\ \hline
\end{tabular}
\caption{Results of DecQN compared with DAPP, RP and GP on small instances}
\label{tab:resultsvrpscd0}
\end{table}

We now discuss the performance of DecQN on instances with Moderate, High and Very High customer density. The test set is obtained by sampling 500 customer realizations and 500 demand realizations for each instance $i \in I$ with Moderate, High and Very High customer density, resulting in 250,000 scenarios.
The aggregated results are reported in Table \ref{tab:resultsvrpscd1}, where the first four columns indicate the characteristics of the instance.
The expected total demand for each instance is reported in column $\mathbb{E}\sum d_c$. The results of DecQN, along with the results of the two benchmarks, are reported in columns DecQN, RP, and GP, respectively. 
\begin{table}[!hbtp]
    \centering
    \footnotesize
    \begin{tabular}{ c c c c c c c c c c c }
        \hline
        $\Psi_{n_z}$ & $Q$ & $m$ & $L$ & $\mathbb{E}\sum d_c$ & DecQN & \% $\mathbb{E}\sum d_c$ & RP  & GP & \%RP & \%GP\\ \hline
        \multirow{3}{*}{Moderate} & 25 & \multirow{3}{*}{3}  & \multirow{3}{*}{221.5} & \multirow{3}{*}{230} & 159.1 & 69.2\% & 99.9  & 143.0 & 59.2\%   & 11.3\% \\
                   & 50 & & & & 193.1 & 83.9\% & 120.6 & 171.4 & 60.2\%   & 12.7\% \\
                   & 75 & & & & 210.1 & 91.4\% & 128.9 & 192.4 & 63.0\%   & 9.2\% \\ \hline
        \multicolumn{9}{r}{\textbf{Avg}} & \textbf{60.8\%}   & \textbf{11.1\%} \\ \hline
        \multirow{3}{*}{High} & 25 & \multirow{3}{*}{7}  & \multirow{3}{*}{195.5} & \multirow{3}{*}{530} & 360.9 & 68.1\% & 217.9 & 321.5 & 65.6\%   & 12.3\% \\
                   & 50 & & & & 452.7 & 85.4\% & 264.5 & 417.9 & 71.1\%   & 8.3\%\\
                   & 75 & & & & 501.2 & 94.6\% & 269.6 & 460.0 & 85.9\%   & 8.9\%\\ \hline
        \multicolumn{9}{r}{\textbf{Avg}} & \textbf{74.2\%}   & \textbf{9.8\%}\\ \hline
        \multirow{3}{*}{Very High} & 25 & \multirow{3}{*}{11} & \multirow{3}{*}{187.3} & \multirow{3}{*}{830} & 552.5 & 66.6\% & 334.0 & 502.5 & 65.4\%   & 9.9\% \\
                   & 50 & & & & 703.0 & 84.7\% & 402.8 & 664.8 & 74.6\% & 5.8\%\\
                   & 75 & & & & 778.7 & 93.8\% & 412.2 & 738.2 & 88.9\%   & 5.5\% \\ \hline
        \multicolumn{9}{r}{\textbf{Avg}} & \textbf{76.3\%}   & \textbf{7.0\%} \\ \hline
        \multicolumn{9}{r}{\textbf{Total Avg}}  & \textbf{70.4\%}   & \textbf{9.3\%} \\ \hline
    \end{tabular}
    \caption{Results of DecQN compared with RP and GP}
    \label{tab:resultsvrpscd1}
\end{table}
DecQN outperforms RP and GR by 70.4\% and 9.3\%, on average. The significant improvement over RP demonstrates the effect of the learning. According to Table \ref{tab:resultsvrpscd1}, \%GP falls from 11.1\% to 7.0\%, as the customer density increases.
As previously stated, we argue that, as the customer density increases, the instance becomes relatively simpler for GP. 
To demonstrate this behavior in our results, we analyze the sensitivity of GP to variations in customer density.
To this purpose, we normalize the demand served by GP by $\mathbb{E}\sum d_c$. Notice that the value of $\mathbb{E}\sum d_c$ solely depends on the level of $\Psi_{n_z}$; hence, $\%\frac{\textrm{GP}}{\sum d_c}$ gives us a normalized measure that can be used to assess the performance of GP with respect to variations in customer density.
This value can also be viewed as the ratio of demands served by GP.
The right side of Table \ref{tab:resultsvrpscd2} reports the average values of $\%\frac{\textrm{GP}}{\sum d_c}$. 
We observe that, as the customer density rises, GP serves a larger ratio of demands. Specifically, the average ratio of demands served by GP increases from 73.4\% to 76.5\% when the customer density increases from Moderate to Very High. 
One intuitive explanation for this is that, as customers in the service area become denser, the problem for each vehicle will be more of a customer selection problem rather than a routing problem, where a long-term collective reward plays a decisive role.
As a result, GP, which ignores  the routing component, performs better on more densely populated service area.
The left side of Table  \ref{tab:resultsvrpscd2} reports average values of $\%\frac{\textrm{DecQN}}{\mathbb{E}\sum d_c}$.
We observe that the performance of the DecQN remains approximately constant as the customer density varies. In particular, the value of $\%\frac{\textrm{DecQN}}{\mathbb{E}\sum d_c}$, averaged across different values of $Q$, remains around $\approx 82\%$ for each level of customer density. This finding is noteworthy because it implies that increases in the problem size due to higher customer density $\Psi_{n_z}$, do not deteriorate the performance of DecQN. Additionally, in line with the earlier discussion, we may attribute this result to the fact that, unlike GP, DecQN accounts for routing.

\begin{table}[!hbtp]
\centering
\footnotesize
\begin{tabular}{cccccc} \hline
\multicolumn{2}{c}{\multirow{2}{*}{$\%\frac{\textrm{DecQN}}{\mathbb{E}\sum d_c}$}} & \multicolumn{3}{c }{$Q$} & \\ \cline{3-5}
\multicolumn{2}{ c }{} & 25 & 50 & 75 & \textbf{Avg} \\ \hline
\multirow{3}{*}{$\Psi_{n_z}$} & Moderate & 69.2\% & 83.9\% & 91.4\% & \textbf{81.5\%} \\ 
 & High & 68.1\% & 85.4\% & 94.6\% & \textbf{82.7\%} \\ 
 & Very High & 66.6\% & 84.7\% & 93.8\% & \textbf{81.7\%} \\ \hline
 & \textbf{Avg} & \textbf{67.9\%} & \textbf{84.7\% }& \textbf{93.3\%} & \textbf{81.9\%} \\ \hline
\end{tabular}
\quad
\begin{tabular}{cccccc}
\hline
\multicolumn{2}{c}{\multirow{2}{*}{$\%\frac{GP}{\mathbb{E}\sum d_c}$}} & \multicolumn{3}{c }{$Q$} & \\ \cline{3-5}
\multicolumn{2}{ c }{} & 25 & 50 & 75 & \textbf{Avg} \\ \hline
\multirow{3}{*}{$\Psi_{n_z}$} & Moderate & 62.2\% & 74.5\% & 83.6\% & \textbf{73.4\%} \\ 
 & High & 60.7\% & 78.9\% & 86.8\% & \textbf{75.4\%} \\ 
 & Very High & 60.6\% & 80.1\% & 88.9\% & \textbf{76.5\%} \\ \hline
 & \textbf{Avg} & \textbf{61.1\%} & \textbf{77.8\%} & \textbf{86.5\%} & \textbf{75.1\%} \\ \hline
\end{tabular}
\caption{The ratio of demands served by DecQN and GP}
\label{tab:resultsvrpscd2}
\end{table}

Finally, Table \ref{tab:resultsvrpscd1} shows that the \%GP is usually higher when $Q$ is smaller. For example, considering instances with Very High customer density, the \%GP declines from 9.9\% to 5.5\% when the vehicle capacity increases from 25 to 75.
A possible reason for this result is that, according to Table \ref{tab:resultsvrpscd2}, as $Q$ increases, the percentage of demands served by DecQN and GP increases to serve almost all demands. 
Thus there is a lower margin of improvement in higher vehicle capacities, which reduces the potential difference between their performances.

\subsection{Performance of DecQN in VRPSD} \label{sec:resultsvrpsd1}
The VRPSD can be seen as a particular version of the \vrp{}, where the set of customers, including their locations and expected demands, is known in advance.
Given that the purpose of DecQN  is to handle variable customer sets, it does not assume the knowledge of customer locations, nor is designed to take advantage of such information. However,  we now show that DecQN can compete with the solution method proposed by \cite{Goodson2016}, which is the best-performing benchmark specialized for VRPSDs with duration limits and multiple vehicles.

In Section \ref{sec:resultsvrpsdinstance}, we define the instance set and describe the scenario generation procedure. For each instance, we train DecQN over 3 million scenarios. Similar to the previous section, we decay the exploration and the learning rates linearly from $1.0$ and $10^{-3}$ to $0.1$ and $10^{-4}$ in the first 300K and 1.5 million trials. We describe the benchmark in Section \ref{sec:resultsvrpsdbench} and discuss the results in Section~\ref{sec:resultsvrpsdresults}.

\subsubsection{VRPSD Instances} \label{sec:resultsvrpsdinstance}
The VRPSD, as a particular version of the \vrp{}, assumes that the number of customers, their locations and their expected demands are given in advance. Therefore, the distribution functions $\Psi_{n_z}, \Psi_l, \Psi_{\bar{d}}$ do not play a role in this section.
The remaining instance-defining components ($\mathcal{A}, \Psi_{\hat{d}}, m, Q, L$) are determined according to the procedure used in \cite{Goodson2016}, which is described in Appendix B.1.2.

\cite{Goodson2016} conduct experiments on 216 instances. However, testing the DecQN on the full set of instances is time-consuming and beyond the scope of this paper. Therefore, we choose a subset of relatively challenging instances. 
In particular, we focus on instances with $n=75$ constructed from R101 Solomon instances. In these instances, the number of vehicles $m$ is 11. We  identify a VRPSD instance $i\in I$ as ($L, Q, U$) with $L\in \mathcal{L}$, $Q\in \mathcal{Q}$, and $U\in\mathcal{U}$. We denote $U$ as the stochastic variability of the demands $\Psi_{\hat{d}}$. We define $\mathcal{L}, \mathcal{Q},$ and $\mathcal{U}$ as follows: $\mathcal{L}=\{\textrm{Short, Medium, Long}\}, \mathcal{Q}=\{25, 50, 75\},$ and $\mathcal{U}=\{\textrm{Low, Moderate, High}\}$. A total of 27 instances are tested. 

\subsubsection{Benchmarks} \label{sec:resultsvrpsdbench}
\cite{Goodson2016} proposed a rollout algorithm~(RA) based on a restocking fixed-route policy. The RA is a form of forward-dynamic programming that can be viewed as a one-step policy iteration \citep{Bertsekas2013}. At each decision epoch, they generate a set of fixed routes using a local search heuristic. Each fixed route is then evaluated using a reward-to-go function for a set of sample demand realizations. Finally, the fixed route with the highest reward-to-go is selected, and the next move on that route is chosen. To account for preemptive actions for instances with up to 25 customers, they solve an auxiliary dynamic program to evaluate the reward-to-go of following a fixed route with and without performing a restocking move. For instances with more than 25 customers, they proposed a dynamic decomposition method to partition customers between vehicles, enabling them to handle instances of up to 100 customers.
\cite{Goodson2016} obtain the initial fixed route by solving the VRPSD using a Simulated Annealing algorithm. They showed that the high-quality initial fixed route plays a vital role. According to their findings, starting with a high-quality initial solution outperforms the case of  starting with a low-quality (randomly produced) fixed route by an average of 11.2\%. We refer to this benchmark policy as GHQ when the algorithm starts with a high-quality fixed route, and refer to it as GLQ when the algorithm starts with a low-quality fixed route.

\subsubsection{Results and Discussion} \label{sec:resultsvrpsdresults}
The test set $\hat{I}(i)$ is obtained by sampling 500 demand realizations for each  instance $i\in I$. The performance of our method is then compared with two benchmark policies, GLQ and GHQ. GLQ is analogous to our method in that, like DecQN, it does not need computing an initial route for each instance. On the other hand, GHQ necessitates the computation of a high-quality initial route for each $i\in I$ instance.

The comparison is summarized in Table \ref{tab:resultsvrpsd1}. In this table, the first three columns indicate the characteristics of the instance. The duration limits S, M, and L, respectively, refer to Short, Medium, and Long. The values L, M, and H for the stochastic variability are abbreviations for Low, Moderate, and High, respectively. The total served demand by DecQN averaged over 500 realizations ($\hat{I}(i)$) for each instance $i\in I$ is reported as DecQN. Column (\%$\mathbb{E}\sum d_c$) shows the ratio of demands that DecQN serves. Note that the total expected demand for each instance is 1079. The results of the benchmarks are shown in columns GLQ and GHQ. Columns \%GLQ and \%GHQ display the percentage improvement of our method over GLQ and GHQ. 
\begin{table}[!hbtp]
\centering
\footnotesize
\begin{tabular}{ccccccccc}
\hline
$L$ & $Q$ & $U$ & DecQN & \% $\mathbb{E}\sum d_c$ & GLQ & GHQ & \%GLQ & \%GHQ \\ \hline
\multirow{9}{*}{S} & \multirow{3}{*}{25} & L & 580.6  & 53.8\%  & 534.6  & 595.9  & 8.6\%  & -2.6\% \\
  &   & M & 556.3  & 51.6\%  & 475.1  & 566.6  & 17.1\% & -1.8\% \\
  &  & H & 534.0  & 49.5\%  & 479.8  & 527.5  & 11.3\% & 1.2\%  \\ \cline{2-9}
  & \multirow{3}{*}{50} & L & 799.7  & 74.1\%  & 738.3  & 834.0  & 8.3\%  & -4.1\% \\
  &  & M & 781.0  & 72.4\%  & 718.5  & 794.7  & 8.7\%  & -1.7\% \\
  &  & H & 751.0  & 69.6\%  & 704.9  & 760.4  & 6.5\%  & -1.2\% \\ \cline{2-9}
  & \multirow{3}{*}{75} & L & 934.0  & 86.6\%  & 848.7 & 1011.4 & 10.1\% & -7.7\% \\
  &  & M & 902.4  & 83.6\%  & 845.2  & 966.2  & 6.8\%  & -6.6\% \\
  &  & H & 885.4  & 82.1\%  & 839.3  & 913.5  & 5.5\%  & -3.1\% \\ \hline
\multirow{9}{*}{M} & \multirow{3}{*}{25} & L & 825.8  & 76.5\%  & 785.6  & 842.9  & 5.1\%  & -2.0\% \\
 &  & M & 800.3  & 74.2\%  & 778.2  & 808.6  & 2.9\%  & -1.0\% \\
 &  & H & 770.7  & 71.4\%  & 747.3  & 776.2  & 3.1\%  & -0.7\% \\ \cline{2-9}
 & \multirow{3}{*}{50} & L & 1054.8 & 97.8\%  & 1019.2 & 1077.5 & 3.5\%  & -2.1\% \\ 
 &  & M & 1035.8 & 96.0\%  & 1002.8 & 1067.0 & 3.3\%  & -2.9\% \\
 &  & H & 1023.8 & 94.9\%  & 977.1  & 1043.6 & 4.8\%  & -1.9\% \\ \cline{2-9}
 & \multirow{3}{*}{75} & L & 1074.5 & 99.6\%  & 1058.3 & 1078.4 & 1.5\%  & -0.4\% \\
 &  & M & 1072.8 & 99.4\%  & 1052.1 & 1077.2 & 2.0\%  & -0.4\% \\
 &  & H & 1069.7 & 99.1\%  & 1044.4 & 1075.3 & 2.4\%  & -0.5\% \\ \hline
\multirow{9}{*}{L} & \multirow{3}{*}{25} & L & 1004.1 & 93.1\%  & 974.4  & 1039.0 & 3.1\%  & -3.4\% \\
 &  & M & 984.5  & 91.2\%  & 959.9  & 998.6  & 2.6\%  & -1.4\% \\
 &  & H & 971.3  & 90.0\%  & 940.3  & 968.0  & 3.3\%  & 0.3\%  \\ \cline{2-9}
 & \multirow{3}{*}{50} & L & 1078.7 & 100.0\%  & 1078.3 & 1078.3 & 0.0\%  & 0.0\%  \\
 &  & M & 1080.3 & 100.1\% & 1076.3 & 1076.9 & 0.4\%  & 0.3\%  \\
 &  & H & 1081.4 & 100.2\% & 1072.3 & 1075.3 & 0.9\%  & 0.6\%  \\ \cline{2-9}
 & \multirow{3}{*}{75} & L & 1079.7 & 100.1\% & 1078.4 & 1078.4 & 0.1\%  & 0.1\%  \\
 &  & M & 1081.8 & 100.3\% & 1077.4 & 1077.5 & 0.4\%  & 0.4\%  \\
 &  & H & 1082.7 & 100.3\% & 1075.5 & 1075.6 & 0.7\%  & 0.7\%  \\ \hline
  &   &   &         &          &         & \textbf{Avg}    & \textbf{4.6\%}  & \textbf{-1.6\%}\\ \hline
\end{tabular}
    \caption{Results of DecQN compared with \cite{Goodson2016}}
    \label{tab:resultsvrpsd1}
\end{table}

\begin{table}[!hbtp]
\centering
\footnotesize
\begin{tabular}{ccccccccccccc} \hline
\multirow{2}{*}{Avg Imp.} & \multicolumn{3}{c}{$L$} & & \multicolumn{3}{c}{$Q$} & & \multicolumn{3}{c}{$U$} & \multirow{2}{*}{Total} \\ \cline{2-4} \cline{6-8} \cline{10-12}
 & S & M & L & & 25 & 50 & 75 & & L & M & H &   \\ \hline
\%GLQ & 9.2\%  & 3.2\%  & 1.4\% & & 6.3\%  & 4.0\%  & 3.3\% & & 4.5\%  & 4.9\%  & 4.3\% & 4.6\% \\
\%GHQ & -3.1\% & -1.4\% & -0.3\% & & -1.3\% & -1.5\% & -1.9\% & & -2.5\% & -1.7\% & -0.5\% & -1.6\% \\ \hline
\end{tabular}

\caption{Average improvement of DecQN over GLQ and GHQ with respect to $L, Q,$ and $U$}
\label{tab:resultsvrpsd11}
\end{table}

Overall, results show that DecQN outperforms the GLQ on average by 4.6\%, while GHQ outperforms DecQN by 1.6\%. To facilitate the discussion, the results are aggregated with respect to $L, Q,$ and $U$ 
in Table \ref{tab:resultsvrpsd11}. 
Comparing the percentage improvements, averaged for each level of duration limit, demonstrates the considerable impact of the duration limit on the gap between our method and two benchmarks.
We observe that,  as the duration limit becomes longer, the gap between the performance of methods is less pronounced.
In particular, when the duration limit changes from Short to Long, the absolute average value of \%GLQ and \%GHQ fall from 9.2\% and 3.1\% to 1.4\% and 0.3\%, respectively.
This reduction in the gap may be explained by observing that, when the duration limit increases, it become less binding, and all methods are able to serve almost all demands.
This logic is also supported by results reported in the fifth column (\%$\mathbb{E}\sum d_c$) of Table \ref{tab:resultsvrpsd1}. The DecQN serves almost all demands (at least 90.0\%, except for instances with $(L,Q)=(\textrm{Medium}, 25)$) when the duration limit is Medium or Long. 
We now asses the impact of the vehicle capacity. Table \ref{tab:resultsvrpsd11} shows that the average \%GHQ declines when vehicle capacity is reduced, whereas DecQN widens the gap with GLQ. Since smaller vehicle capacities imply a larger number of restocking operations, this result may be interpreted as a sign that DecQN successfully handles preventative restocking operations.

\begin{figure}
    \centering
    \includegraphics[width=0.6\textwidth]{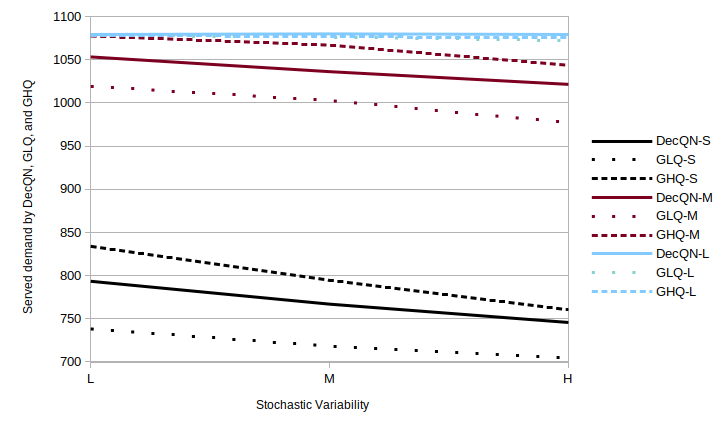}
    \caption{Performance of DecQN, GLQ, and GHQ with respect to the stochastic variability}
    \label{fig:resultsvrpsd}
\end{figure}

Regarding the impact of the stochastic variability, Table \ref{tab:resultsvrpsd11} shows an overall reduction in the absolute average gaps (\%GLQ and \%GHQ) as $U$ increases from Low to High. This behavior can be explained by looking into the difference between GLQ and GHQ. 
Figure \ref{fig:resultsvrpsd} illustrates the trend of served demands by DecQN, GLQ, and GHQ, when the capacity is 50, in different levels of duration limit with respect to changes in stochastic variability. 
In this figure, results for different duration limits are identified by an indication appended to the policy name.
For example, DecQN-S refers to the results of DecQN in instances with a Short duration limit. As shown in this figure, the increase in $U$ narrows the gap between GLQ and GHQ when the duration limit is Low. 
This finding implies that as the stochastic variability increases, having a high-quality initial solution becomes less important. This figure also explains why the gap between GLQ and GHQ, and consequently \%GLQ and \%GHQ, do not follow a descending trend (as in instances with Short $L$) when $U$ rises. 
The amount of served demands in instances with Medium and Long duration limit are very close to the total expected demand (i.e., 1079), making the performance of policies bounded.
Hence the gap between GLQ and GHQ cannot follow the same trend as in instances with a Short duration limit.
Another important finding is that while the performance of GHQ, when $U$ increases from Low to High in instances with Short $L$, drops by 10.0\% (from 595.9, 834.0, and 1011.4 to 527.5, 760.4, and 913.5, respectively for Small, Medium, and Large $Q$), the performance of our DecQN decreases only by 6.4\% (from 580.6, 799.7, and 934.0 to 534.0, 751.0, and 885.4, respectively for Small, Medium, and Large $Q$). It can be deduced that our method is relatively better at dealing with stochastic variability than GHQ.

A last comment concerns the computing times. Given that GLQ starts from a random set of routes, the initialization time is negligible. Conversely, GHQ necessitates the computation of a high-quality initial route for each instance $i \in I$.
The computing times for establishing initial routes are not reported in \cite{Goodson2016}. However, results shared with us by the authors, we conclude that the average CPU time dedicated to obtaining  high-quality initial solutions for instances with 75 customers was 25,466 seconds ($\approx 7$ hours). In comparison, our method,  once trained for 24 to 36 hours, can be used with no further computation as long as the instance (including the associated distribution functions and parameters such as $L, Q, U$) remains the same.
Therefore, we can roughly claim that using our method for more than five ($\approx \frac{36}{7}$) operating days is computationally advantageous over the GHQ.

Summarizing, we observed that DecQN and GLQ are analogous in the sense they both start with a random initial policy and solution. Results showed that the DecQN method could outperform GLQ by 4.6\% on average. On the other hand, the average gap between DecQN and GHQ, where a high-quality initial solution is required, is -1.6\%. However, results showed that the advantages of starting from a high-quality initial solution become less pronounced for higher values of stochastic variability. Furthermore, obtaining high-quality initial solutions is computationally expensive, and DecQN resulted computationally advantageous with respect to GHQ after five days of operations.

\subsection{Generalized DecQN} \label{sec:resultsaggr1}
The previous section showed that when DecQN is trained on VRPSD instances, it can compete with \cite{Goodson2016}. However, it was necessary to train a new policy for each instance $i\in I$.
In this experiment, we further explore the flexibility of the DecQN in terms of its generalization capabilities.  Specifically, we investigate the possibility to train one DecQN model to solve a group of VRPSD instances $I'$, where $I'\subset I$. This also has some practical advantages, as it would be beneficial for a manager to apply the same policy regardless of, for example, the duration limit.

We use the same instances described in Section \ref{sec:resultsvrpsdinstance}, but we restrict our analysis to instances with $Q=50$ (i.e., $I=\mathcal{L}\times \{50\} \times \mathcal{U}$).
We carry out the experiment in three steps. In the first step, we aggregate instances with respect to the duration limit $L$. Accordingly, we train three policies, each for a $U\in\mathcal{U}$. More specifically, in instances aggregated over $L$, the duration limit follows a uniform distribution function to pick a value from $\mathcal{L}$. Similarly, we generalize the method over the stochastic variability $U$ in the next step. In the last step, we generalize the DecQN over both components. In this setting, a single policy is trained to solve all nine instances. To provide enough information to DecQN, we append the value of the component on which we are aggregating to the observation representation. For example, when aggregating over the duration limit, we pass $L$ as an extra component into Equation~\eqref{eq:observation}, resulting in $\displaystyle{o_{k,\bar{v}}=[F_k|H_k|G_k|t_k|L]}$.

Once we obtain a policy, we test it on all instances considered in the training phase. As in the previous experiments, the training and testing realization sets are disjoint. For the first step, Table~\ref{tab:aggrdl} summarizes the results of the trained policies ($\textrm{DecQN}_{L}$) and compares them with those trained in Section \ref{sec:resultsvrpsd1} (DecQN) and the two benchmark policies, GLQ and GHQ. Detailed computational results are reported in Table 12, in Appendix C. Table \ref{tab:aggrdl} reports the improvement percentage of the generalized policy over the benchmarks, averaged for each instance characteristic.
Results show that the aggregation over the duration limit performs on average 2.3\% and 3.7\% worse than the individual policies (DecQN) and the GHQ, respectively. However, it still outperforms the GLQ by 1.6\%. It can also be seen that the gap between $\textrm{DecQN}_{L}$ and DecQN decreases when the duration limit gets longer (i.e., from -4.6\% to -0.1\%). Also, the aggregation works relatively better in higher stochastic variability.
\begin{table}[!hbtp]
\centering
\footnotesize
\begin{tabular}{cccccccccc} \hline
\multirow{2}{*}{Avg Imp.} & \multicolumn{3}{c}{L} & & \multicolumn{3}{c}{U} & \multirow{2}{*}{Total} \\ \cline{2-4} \cline{6-8}
 & S & M & L & & L & M & H & \\ \hline
\%DecQN & -4.6\% & -2.2\% & -0.1\% & & -2.3\% & -2.5\% & -2.0\% & -2.3\%                \\
\%GLQ   & 2.9\%  & 1.6\%  & 0.3\%  & & 1.5\%  & 1.4\%  & 1.9\%  & 1.6\%                 \\
\%GHQ   & -6.8\% & -4.4\% & 0.2\%  & & -4.3\% & -3.9\% & -2.8\% & -3.7\%    \\ \hline           
\end{tabular}
\caption{Average improvement of $\textrm{DecQN}_{L}$ over DecQN, GLQ, and GHQ}
\label{tab:aggrdl}
\end{table}

Table \ref{tab:resultsaggrsv} reports the results of the second step ($\textrm{DecQN}_{U}$), where instances are aggregated with respect to the stochastic variability. Detailed  results are provided in Table 13 in Appendix C. Similar to the previous analysis, we compare the results with the individual (DecQN) and benchmark (GLQ and GHQ) policies. Results show that our method can handle a broader range of instances at the expense of performing on average 0.7\% worse than the individual policies. Interestingly, it can be seen that the generalized DecQN performs as well as the individual policies in instances with Medium and Long duration limits. Although the generalized model has a gap of -2.2\% with the GHQ, it still considerably outperforms the GLQ (3.2\%).

\begin{table}[!hbtp]
\centering
\footnotesize
\begin{tabular}{cccccccccc} \hline
\multirow{2}{*}{Avg Imp.} & \multicolumn{3}{c}{L} & & \multicolumn{3}{c}{U} & \multirow{2}{*}{Total} \\ \cline{2-4} \cline{6-8}
 & S & M & L & & L & M & H & \\ \hline
\%DecQN & -2.2\% & 0.0\%  & 0.0\% & & -1.2\% & -0.8\% & -0.3\% & -0.7\% \\
\%GLQ & 5.4\%  & 3.9\%  & 0.4\% & & 2.7\%  & 3.3\%  & 3.8\%  & 3.2\% \\
\%GHQ & -4.5\% & -2.3\% & 0.3\% & & -3.2\% & -2.2\% & -1.1\% & -2.2\% \\ \hline
\end{tabular}
\caption{Average improvement of $\textrm{DecQN}_{U}$ over DecQN, GLQ, and GHQ}
\label{tab:resultsaggrsv}
\end{table}

Finally, Table \ref{tab:aggrall} reports the performance of the single policy trained on all nine instances ($\textrm{DecQN}_{All}$). Detailed results are provided in Table 14 in Appendix C. Similar to Tables \ref{tab:aggrdl} and \ref{tab:resultsaggrsv}, we compare the performance of the generalized DecQN with the individual and benchmark policies. Results show that $\textrm{DecQN}_{All}$ is able to solve a set of instances with varying values of $L$ and $U$ at the expense of a 2.1\% reduction of the performance quality.
Confronting these results with those obtained in Tables \ref{tab:aggrdl} and \ref{tab:resultsaggrsv}, we can deduce that most of the performance deterioration comes from the generalization over the time limit dimension $L$. 
This experiment also shows that, if suitably trained, the proposed algorithm can be adapted to handle problems with varying features such as, for example, different duration limits.

\begin{table}[!hbtp]
\centering
\footnotesize
\begin{tabular}{cccccccccc} \hline
\multirow{2}{*}{Avg Imp.} & \multicolumn{3}{c}{L} & & \multicolumn{3}{c}{U} & \multirow{2}{*}{Total} \\ \cline{2-4} \cline{6-8}
 & S & M & L & & L & M & H & \\ \hline
\%DecQN & -4.3\% & -1.9\% & -0.2\% & & -2.7\% & -2.2\% & -1.4\% & -2.1\% \\
\%GLQ   & 3.2\%  & 1.9\%  & 0.2\%  & & 1.0\%  & 1.8\%  & 2.5\%  & 1.8\% \\
\%GHQ   & -6.6\% & -4.2\% & 0.1\% & & -4.7\% & -3.6\% & -2.3\% & -3.5\% \\ \hline
\end{tabular}
\caption{Average improvement of $\textrm{DecQN}_{All}$ over DecQN, GLQ, and GHQ}
\label{tab:aggrall}
\end{table}

\section{Conclusions}\label{conclusion}
Motivated by a distribution planning problem arising in domestic donor collection services, we introduced the \vrp{}.
In this problem, customer locations and demands are stochastic, and a fleet of capacitated vehicles must complete the service within a specified time limit. The goal is to maximize the total served demand. We provided two problem formulations based on Markov Decision Process. The first follows a traditional centralized decision-making perspective, while the second models a partially decentralized decision-making framework where vehicles autonomously establish their routes according to the observed state of the system. This second formulation enabled us to drastically reduce the dimension of the state and the action spaces, thus resulting in a more tractable problem.
Variable-sized stochastic customer sets have been addressed via a heatmap-style representation of the customer demands. 
This approach, in particular, partitions the stochastic customers geographically and encodes them as a fixed-size vector. We then solved the resulting problem by a Q-learning algorithm. For this purpose,  we developed a two-layer artificial neural network to approximate state-action Q factors. 
Computational results showed that our method significantly outperforms three benchmarks policies (DAPP, RP, and GP). Although DAPP provided better bounds than the other benchmarks, it could only handle instances with up to 15 customers.
Furthermore, we evaluated the DecQN on the VRPSD, which is a particular case of \vrp{}, and showed that our method could compete with the state-of-the-art solution method specialized for that problem. 
As a further investigation into the capabilities of our framework, we tested and showed that the DecQN could be generalized over several problem components, such as the stochastic variability of demands. Finally, as shown in Appendix D, the obtained policies and value functions can be easily employed as a base policy and a reward-to-go estimator in on-line rollout methods.
 
As DecQN can efficiently handle problems with stochastic customer sets, it is promising to investigate how it could be extended to tackle dynamic versions of the VRP, where customers arise dynamically during the execution period.
Another potential avenue is to investigate how the partially decentralized approach can be applied to other SVRPs, to handle, for example, stochastic travel time.

\end{document}


\section{Algorithms}\label{app:algs}
\subsection{Q-learning algorithm} \label{app:qlearning}
Q-learning is a model-free algorithm based on approximate value iteration in which an agent learns how to behave in order to maximize the collected reward through a sequence of interactions with the environment, resulting in a policy that satisfies Equation (26). More specifically, this method iteratively updates the value of making the action $x_k$ in state $s_k$, known as the Q factor of that state-action pair ($Q(s_k,x_k)$). In fact, it is possible to rewrite the Bellman equation (defined in Equation (25)) in terms of Q factors as follows:
\begin{equation} \label{eq:vqlink} \addtocounter{equation}{31}
    V(s_k) = \max_{x_k\in A(s_k)} Q(s_k, x_k),
\end{equation}
where:
\begin{equation*}
    Q(s_k, x_k) = R(s_k, x_k, w) + \gamma\max_{x_{k+1}\in A(s_{k+1})} Q(s_{k+1}, x_{k+1}).
\end{equation*} 
The reader is referred to \citep{Powell2011} for a complete overview of the algorithm.
Algorithm \ref{alg:qlearning} illustrates a generic Q-learning framework. To summarize, this algorithm simulates the problem for several trials ($Trials\_MAX$) over different sample scenarios of stochastic components of the problem. The simulation follows the so-called $\epsilon$-greedy policy, which enables random exploration of the action space. The level of exploration is typically larger in the initial phases of the algorithm when actions are mostly chosen randomly. As the algorithm proceeds, actions with higher estimated values are generally preferred.
Q-learning exploits the interactions between the agent and the environment to update the Q values (the so-called \textit{experiences}). An experience is defined as a tuple of $(s_k, x_k, r_k, s_{k+1})$, describing interaction of the agent with the environment consisting of the agent observation of the state $s_k$, its action $x_k$, the received reward $r_k$ (in \vrp{}, $r_k=R(s_k, x_k, w)$) and the observation of the next state $s_{k+1}$. For each state-action pair experienced, we estimate their Q value as $\hat{Q}(s_k, x_k)$ using the following equation:
\begin{equation}
    \hat{Q}(s_k, x_k) = r_k + \gamma\max_{x_{k+1}\in A(s_{k+1})} Q(s_{k+1}, x_{k+1}),
    \label{eq:qestimate}
\end{equation}
where $\gamma$ is the discount factor. Then, we update the $Q(s_k, x_k)$ according to the following equation:
\begin{equation}
    Q(s_k, x_k) \leftarrow Q(s_k, x_k) + \alpha * [\hat{Q}(s_k, x_k) -  Q(s_k, x_k)],
    \label{eq:qvalue}
\end{equation}
where $\alpha$ is the step size (learning rate).
\begin{algorithm}[!hbtp]
\footnotesize
 \SetAlgoLined
 \DontPrintSemicolon
 \LinesNumbered
 Initialize: set Q values to 0\;
 \While{trials $<$ Trials\_MAX}
 {
 Choose a sample scenario $w$\;
 Observe state $s_0$\;
 $k\gets 0$, $simulate\gets$True\;
 \While{simulate}{
 Take action: $x_k\gets\begin{cases}
     \textrm{random action in }A(s_k) & \epsilon \\
     \arg\max_{x'\in A(s_k)}Q(s_k, x') & 1-\epsilon
     \end{cases}$\; \vspace{5pt}
 Transit to $s_{k+1} \gets S^M(s_k, x_k, w_{k+1})$ and observe the reward $r_k$ \;
  \vspace{5pt}
 Estimate and update new Q values:\;
 \begin{itemize}
     \item $\hat{Q}(s_k, x_k) \gets \begin{cases}
     r_k & A(s_{k+1})=\emptyset \\
     r_k + \gamma \max_{x'\in A(s_{k+1})} Q(s_{k+1}, x') & \textrm{otherwise} 
     \end{cases}$\; \vspace{5pt}
     \item $Q(s_k, x_k)\gets Q(s_k, x_k) + \alpha [\hat{Q}(s_k, x_k) - Q(s_k, x_k)]$\;
 \end{itemize}
  \vspace{5pt}
 $k\gets k+1$\;
 \If{$A(s_k)=\emptyset$}{
 $simulate\gets$False\;
 }
 }
 Set $trials = trials + 1$\;
 Decay the exploration rate $\epsilon$\;
 }
 \textbf{return} Q values
 \caption{Generic Q-learning Algorithm}
 \label{alg:qlearning}
\end{algorithm}

Algorithm \ref{alg:qlearning} returns a set of Q values that can be used as the decision policy. In particular, the optimal decision policy will be:
\begin{equation} \label{eq:qpolicy}
    \pi^*(s_k)=\arg\max_{x_k\in A(s_k)}Q(s_k, x_k).
\end{equation}

\subsection{The DecQN}\label{app:decqn}
This section illustrates the pseudo code for our proposed DecQN algorithm.
\begin{algorithm}[!hbtp]
\footnotesize
 \SetAlgoLined
 \DontPrintSemicolon
 \LinesNumbered
 Initialize the Q-network and the target Q-network with random weights $\theta, \bar{\theta}$ \;
 $B\gets[]$\;
 \While{trials $<$ Trials\_MAX}
 {
 Choose a sample scenario $w$ and observe a new set of customers $\mathcal{C}_0$\;
 $t_0\gets0$, $k\gets0$, $\bar{v}_{prev}\gets\emptyset$, $\textrm{Buff}\gets[\emptyset]_{v\in\mathcal{V}}$, and $simulate\gets True$\;
 
 \While{simulate}{
 $\bar{v}\gets$ a random $v\in\bar{\mathcal{V}}$, where $\bar{\mathcal{V}} = \{v\in\mathcal{V}|a_v=t_k \}$\;
 \If{$k > 0$}{
 Transit from $s_{k-1}^x$ to $s_k$ 
 and observe the reward $r$ \;
 $\textrm{Buff}[\bar{v}_{prev}] \gets (s_{k-1},\bar{v}_{prev}, x_{k-1}, s_{k}, \bar{v})$\;
 \If{\textup{Buff} $[\bar{v}]\neq \emptyset$}{
 $(s, v, x, s', v')\xleftarrow{} \textrm{Buff}[\bar{v}]$\;
 Append $(s,v, x, r, s', v')$ to $B$\;
 }
 }
 Observe $o_{k, \bar{v}} = O(s_k, \bar{v})$ \;
 Take action: 
 $x_k\gets\begin{cases}
     \textrm{random action in }A(s_k, \bar{v}) & rand[0,1] < \epsilon \\
     \arg\max_{x'}Q(o_{k, \bar{v}}, x', \theta) & \textrm{otherwise}
     \end{cases}$\;
 Transit from $s_k$ to $s_k^x$ \;
 \If{$rand[0,1] < \beta_t$}{
 Estimate and update new Q values:\;
 \begin{itemize}
     \item Sample a batch of experiences $\Tilde{B}$ from $B$
     \item Compute the estimation $\hat{Q}(O(s, v), x, \theta)\gets\begin{cases}
     r & A(s_k, v)=\{l_0\} \land l_v = l_0 \\
     r + \gamma \max_{x'}Q(O(s', v'), x', \bar{\theta}) & \textrm{otherwise}
     \end{cases},~\forall (s, v, x, r, s', v')\in \Tilde{B}$
 \item Compute the loss function: $\Phi\gets\mathbb{E}_{\Tilde{B}}\left[Huber\left(Q\left(O(s, v), x, \theta\right) - \hat{Q}\left(O(s, v), x, \theta\right)\right)\right]$\;
 \item Update the network weights $\theta$:
 $\theta \gets \theta - \alpha \nabla_{\theta} \Phi$\;
 \end{itemize}
 }
 $k \gets k + 1$, and $\bar{v}_{prev} \gets \bar{v}$\;
 \If{$A(s_k, \bar{v})=\emptyset \land l_v=l_0, \forall v\in\mathcal{V}$}{
 $simulate\gets$ False\;
 }
 }
 $trials \gets trials + 1$\;
 Every $\beta_d$ $trials$, $\bar{\theta} \gets \theta$\;
 Decay the learning rate $\alpha$ and the exploration rate $\epsilon$\;
 }
 \textbf{return} Q values\;
 \caption{The DecQN}
 \label{alg:DecQN}
\end{algorithm}

\begin{figure}
    \centering
    \addtocounter{figure}{2}
    \includegraphics[width=0.8\textwidth]{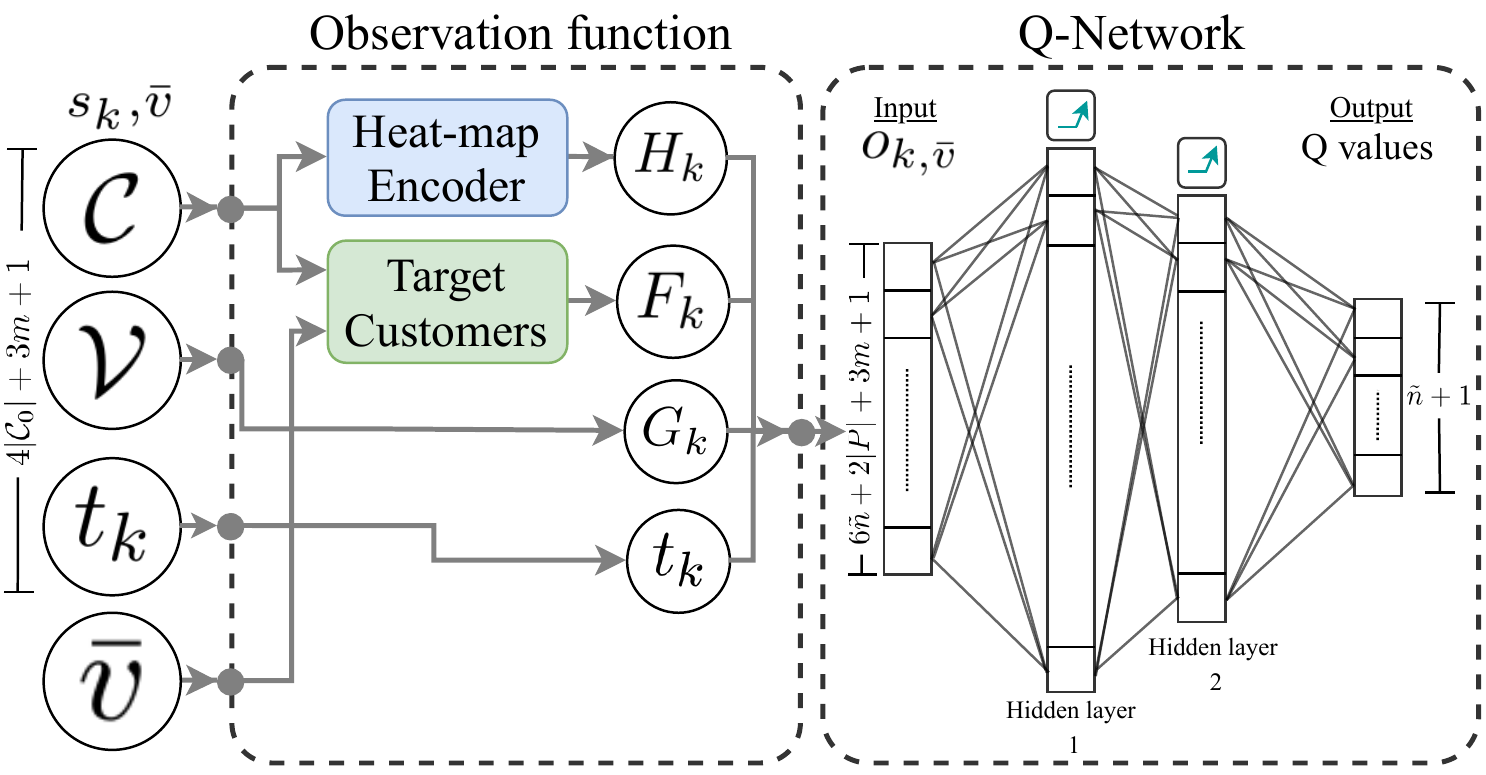}
    \caption{A schematic overview of DecQN}
    \label{fig:overview}
\end{figure}

\newpage
\section{Experimental setup}
\subsection{Instance generation}
This section provides details on the methods used to generate the \vrp{} and VRPSD instances in our experiments.
\subsubsection{VRP-VCSD instances} \label{app:vrpscdinstance}
\begin{figure}[!hbtp]
\centering
\includegraphics[width=5cm]{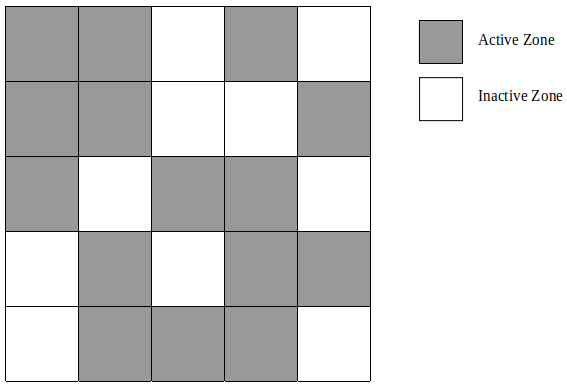}
\caption{The layout of active zones}
\label{fig:activezones}
\end{figure}

We consider a squared-shape service area $\mathcal{A}$ of dimension $100\times 100$ with a depot located at the center ($l_0=(50, 50)$). Customers can potentially be located at any position inside the service area. However, in real-life applications, certain zones of the service area never generate requests, such as non-residential areas.
Therefore, for a given service area partitioned into a set of $5\times 5$ square-shaped zones $Z$, we randomly pick a subset $\mathcal{Z}$ of \emph{active} zones, in which customers must be located. In our experiment, we consider 15 active zones out of 25 (i.e., 60\% of the zones are active). 
Figure \ref{fig:activezones} shows the layout of the considered region.

To generate the customer locations and demands, we make the assumption that active zones are homogeneous in the sense that the customer density $\Psi_{n_z}$ and distributions $\Psi_l, ~\Psi_{\bar{d}},~ \Psi_{\hat{d}}$ do not depend on the specific zone. The customer density $\Psi_{n_z}$ specifies the probability distribution of the number $n_z$ of requests for each active zone.
We consider five possible levels of customer density $\mathcal{D}=\{\textrm{Very Low, Low, Moderate, High, Very High}\}$. 
In particular, $n_z$ takes values \{0, 1, 2\} for Very Low and Low with probabilities \{$0.5, 0.\Bar{3}, 0.1\bar{6}$\} and \{$0.\Bar{3}, 0.\Bar{3}, 0.\Bar{3}$\}, respectively. For denser instances, it takes values \{0, 1, 2, 3\} for Moderate, \{2, 3, 4, 5\} for High, and \{4, 5, 6, 7\} for Very High customer densities, with probabilities \{0.1, 0.4, 0.4, 0.1\}, respectively. 
The corresponding expected number of customers at each zone $\bar{n}_z$ is $0.\bar{6}$, 1.0, 1.5, 3.5, and 5.5, respectively. 
The total expected number of customers $\bar{n}$ for a given $\Psi_{n_z}$ is calculated as $\bar{n}=\mathbb{E} [\sum_{z\in \mathcal{Z}} n_z] = |\mathcal{Z}|*\bar{n}_z$, resulting in 10, 15, 23, 53, and 83 customers for Very Low, Low, Moderate, High, and Very High customer density, respectively.
Once $n_z$ is sampled according to the distribution $\Psi_{n_z}$ for each zone $z \in \mathcal{Z}$, the customer locations are generated following the distribution $\Psi_l$. In this experiment, we choose $\Psi_l$ to be the uniform distribution function over (x,y)-coordinates inside a given zone.
The above process allows us to sample the initial set of customers $\mathcal{C}_0$.

Given a customer $c \in \mathcal{C}_0$, the expected demand $\bar{d}_c$ and actual demand $\hat{d}_c$ are sampled according to distributions $\Psi_{\bar{d}}$ and $\Psi_{\hat{d}}$, respectively. These distributions are chosen following the method proposed in \cite{Gendreau1995}. In particular, $\Psi_{\bar{d}}$ is the discrete uniform distribution over the values $\{5, 10, 15\}$. Once the value $\bar{d}_c$ has been sampled, the actual demands follow the distribution $\Psi_{\hat{d}}$ which is the discrete uniform distribution on values $\{\bar{d}_c-5, ..., \bar{d}_c, ..., \bar{d}_c + 5\}$.
In order to avoid the no-demands case when $\bar{d}_c=5$, we use the discrete uniform distribution function on values $\{\bar{d}_c-4, ..., \bar{d}_c, ..., \bar{d}_c + 4\}$.

We determine the number of vehicles $m$ and the time limit $L$ uniquely as functions of the distribution $\Psi_{n_z}$ by linking them through the concept of the filling rate.
The filling rate $f=\frac{\mathbb{E}[\sum_{c\in \mathcal{C}_0} \bar{d}_c]}{m*Q}$ indicates the portion of all demands that can be served by all vehicles without replenishment. 
In the equation above, $\mathbb{E}[\sum_{c\in \mathcal{C}_0} \bar{d}_c]$ is the expected total demand and can be approximated by $\Bar{n}_z*|\mathcal{Z}|*\bar{\bar{d}}$, where $\bar{\bar{d}}$ is the average expected demand of a customer. Since the distribution function $\Psi_{\bar{d}}$ assigns expected demands 5, 10, or 15 to customers uniformly (with a probability of $\frac{1}{3}$), the average expected demand of a customer ($\bar{\bar{d}}$) will be equal to $10$ ($\frac{5+10+15}{3}$) in all instances. By fixing the filling rate and the vehicle capacity to 1 and 75, respectively, the number of vehicles for each level of customer density $\Psi_{n_z}$ can then be obtained by the definition of the filling rate:
\begin{equation}
    m=\ceil*{\frac{\mathbb{E}[\sum_{c\in \mathcal{C}_0} \bar{d}_c]}{f*Q}}=\ceil*{\frac{\Bar{n}_z*|\mathcal{Z}|*\bar{\bar{d}}}{f*Q}}=\ceil*{\frac{\Bar{n}_z*15*10}{1*75}}=\ceil*{2*\Bar{n}_z}.
    \label{eq:vrpscdm}
\end{equation}
In this way, the number of required vehicles for each class of customer density is 2, 2, 3, 7, and 11, respectively. Notice that, according to the expression $\mathbb{E}[\sum_{c\in \mathcal{C}_0} \bar{d}_c]$, the expected total demand at each instance will be 100, 150, 230, 530, and 830 for Very Low, Low, Moderate, High and Very High customer density, respectively.

To determine the duration limit $L$ for a given instance with the customer density distribution $\Psi_{n_z}$ and its associated $m$, we use the average travel time of vehicles acquired by solving a set of 250,000 realizations of that instance via GP (described in Section 6.1.2). In this process, we assumed that $L=\infty,~Q=75$, and realizations are randomly sampled from distributions $\Psi_{n_z}, \Psi_l, \Psi_{\bar{d}}, \Psi_{\hat{d}}$. To make the duration limit constraining, we let the duration limit be the average travel time of all vehicles multiplied by 0.75. 
The resulting duration limits for Very Low, Low, Moderate, High and Very High customer densities are 143.71, 201.38, 221.47, 195.54, and 187.29, respectively.
To complete the description of the instance generation, we need to determine the values $Q$. 
According to Equation \eqref{eq:vrpscdm}, the number of vehicles $m$ is determined such that when $Q=75$, vehicles can serve the total expected demand without replenishing at the depot. Therefore, for a given duration limit, while using smaller capacities ($<75$) may increase the frequency of restocking operations, which results in serving fewer customers in a limited time, considering larger vehicle capacities ($>75$) may not necessarily affect the performance of vehicles. Therefore, in this experiment, we consider three possible values $\mathcal{Q}=\{25, 50, 75\}$.

To summarize, the instance set $I$ is obtained by varying the customer density distribution function $\Psi_{n_z} \in \mathcal{D}= \{\textrm{Very Low, Low, Moderate, High, Very High}\}$ and the vehicle capacity $Q \in \mathcal{Q} = \{25, 50, 75\}$, resulting in 15 different instances. Fixing the customer density distribution also uniquely determines the number of vehicles $m$ and the time limit $L$. The rest of the instance-defining elements, namely the service area $\mathcal{A}$ and the distributions $\Psi_l, \Psi_{\bar{d}}, \Psi_{\hat{d}}$ are kept fixed throughout the experiment. In other words an instance in this experiment is completely identified  by the pair ($\Psi_{n_z}, Q$). Finally, realizations $\hat{i}$ of an instance $i$ are obtained by sampling from distributions $\Psi_{n_z}, \Psi_l, \Psi_{\bar{d}}, \Psi_{\hat{d}} $.

\subsubsection{VRPSD instances} \label{app:vrpsdinstance}
\cite{Goodson2013} used the \textit{R101} and \textit{C101} instances from \cite{Solomon1987} and took the first $n$ customers to generate the set $\mathcal{C}_0$ with $n=$ 25, 50, 75, and 100.
Concerning expected demands, the authors used the deterministic demands provided by the Solomon's instances. 
To determine the number of vehicles $m$ and the duration limit $L$, the authors solved each instance as a classical VRP using a heuristic with a vehicle capacity of 100. The resulting number of routes determines $m$, while a set of three time limits are obtained by multiplying the length of the longest route by 0.75, 1.25, and 1.75, resulting in  $\mathcal{L}=\{\textrm{Short, Medium, Long} \}$. The  considered duration limits Short, Medium, and Long are 103.05, 171.75, and 240.45, respectively. They also varied the vehicle capacity $Q$ to take a size in $\mathcal{Q}=\{25, 50, 75 \}$. 

Regarding the distribution functions of the customer demands, they characterized $\Psi_{\hat{d}}$ by an instance identifier $U$, denoted as the stochastic variability of the demands. The stochastic variability $U$ may have a value in $\mathcal{U}=\{\textrm{Low, Moderate, High} \}$. In particular, $\Psi_{\hat{d}}$ takes the distribution function of $\{\bar{d}/2, \bar{d}, 3\bar{d}/2 \}$ with probabilities $(0.05, 0.9, 0.05)$ for Low, $\{0, \bar{d}/2, \bar{d}, 3\bar{d}/2, 2\bar{d} \}$ with probabilities $(0.05, 0.15, 0.6, 0.15, 0.05)$ for Moderate, and $\{0, \bar{d}/2, \bar{d}, 3\bar{d}/2, 2\bar{d} \}$ with a uniform probability for High stochastic variability, respectively. For details, the reader may refer to \cite{Goodson2016}. In consequence, the authors conduct experiments on 216 instances, comprised of eight different sets of customers (i.e., four levels for the number of customers $\{25, 50, 75, 100 \}$ $\times$ two classes of instances $\{\textrm{R101, C101}\}$). In particular, for each set of customers with the same locations and expected demands, they did tests on 27 distinct instances.

\subsection{Hyper-parameters and computing facility} \label{app:hyperparameters}
The hyper-parameters in Algorithm~\ref{alg:DecQN} have been set as follows. We consider a minibatch of 32 experiences ($|\Tilde{B}|=32$) uniformly sampled from the replay memory with a size of 50000, where $\beta_t=0.05$. We also replace the target network every 1000 trials ($\beta_d=1000$). For all experiments, we set the discount factor $\gamma$ to 0.999~\citep{Kullman2020}. The hyper-parameters related to the maximum number of trials ($Trials\_MAX$) and the decaying factor for the learning rate ($\alpha$) and the exploration rate ($\epsilon$) will be specified for each experiment later. All the procedures have been implemented in python and executed on 3.6 GHz Intel Xeon processors with 64 GB RAM (no GPUs are used).

\subsection{Deterministic Formulation}\label{app:mf}
In this section, we present a mathematical formulation for the deterministic version of the \vrp{}. Unlike DecQN, where both customers and demands are stochastic, this formulation assumes that the set of customers is known, and the actual demand is equal to the expected value. In order to provide an upper bound for the possible number of trips for a vehicle ($E$), we propose a heuristic algorithm, which is illustrated in Algorithm \ref{alg:cap}. This algorithm computes the total number of trips for a vehicle in the worst-case scenario when it visits customers one-per-trip from closest to farthest in terms of their distance to the depot. This strategy results in the highest number of trips.
\begin{table}[!h]
    \centering
    \addtocounter{table}{9}
    \footnotesize
    \begin{tabular}{c|l}
        Notation & Description \\ \hline
        $N$ & set of nodes, $N=\mathcal{C}_0\cup \{0\}$. Node 0 represents the depot. \\
        $i, j$ & index of nodes in $N$\\
        $\mathcal{E}$ & set of trips for a vehicle, $e\in \mathcal{E}$, $\mathcal{E}=\{0, 1, ..., E-1\}$ \\
        $E$ & an upper bound for the number of trips for a vehicle\\
        $M$ & a big number\\
        \hline
        $x^{ve}_{i}$ & the amount of node $i$'s demand that is served by vehicle $v$ in trip $e$ \\
        $y^{ve}_{ij}$ & binary variable takes 1 if node $j$ is immediately visited after node $i$ by vehicle $v$ in trip $e$ \\
        $\lambda^{ve}_{i}$ & binary variable takes 1 if customer $i$ is visited by vehicle $v$ in trip $e$ \\
        $q^{ve}_{i}$ & available capacity of vehicle $v$ when it arrives at node $i$ in trip $e$ \\
        $t^{ve}_{i}$ & time at which the vehicle $v$ arrives at node $i$ in trip $e$ \\
        $\underline{t}^{ve}$ & time at which the vehicle $v$ leaves the depot for trip $e$ 
        
    \end{tabular}
    \caption{Notations of the deterministic formulation}
    \label{tab:my_label}
\end{table}
\begin{align}
    &\max \qquad && \sum_{i\in N}\sum_{v\in \mathcal{V}}\sum_{e\in \mathcal{E}} x^{ve}_{i} && \label{eq:mfc0}\\
    & \text{s.t.} \qquad && \sum_{j\in N} y^{ve}_{ij} \leq 1; && \forall i\in N, v\in \mathcal{V}, e\in \mathcal{E} \label{eq:mfc1}\\
    &&&\sum_{j\in N} y^{ve}_{ji} \leq 1; && \forall i\in N, v\in \mathcal{V}, e\in \mathcal{E} \label{eq:mfc2}\\
    &&&y^{ve}_{ii} = 0;&& ~\forall i\in N, v\in \mathcal{V}, e\in \mathcal{E} \label{eq:mfc3}\\
    &&&\sum_{i\in N} y^{ve}_{0i} \leq \sum_{i\in N} \lambda^{ve}_{i};&& ~\forall v\in \mathcal{V}, e\in \mathcal{E} \label{eq:mfc4}\\
    &&&\sum_{j\in N} y^{ve}_{ji} = \sum_{j\in N} y^{ve}_{ij} = \lambda^{ve}_{i};&& ~\forall i\in N, v\in \mathcal{V}, e\in \mathcal{E} \label{eq:mfc5}\\
    &&&\sum_{i\in N} \sum_{j\in N} y^{ve}_{ij} \leq M(\sum_{i\in N} \sum_{j\in N} y^{v, e-1}_{ij});&& ~\forall v\in \mathcal{V}, e\in \mathcal{E}\setminus 0 \label{eq:mfc6} \\
    &&&\sum_{i\in N} \lambda^{ve}_{i} \leq M(\sum_{i\in N} \lambda^{v, e-1}_{i});&& ~\forall v\in \mathcal{V}, e\in \mathcal{E}\setminus 0 \label{eq:mfc7}\\
    &&&x^{ve}_{i} \leq M \lambda^{ve}_{i};&& ~\forall i\in N, v\in \mathcal{V}, e\in \mathcal{E}  \label{eq:mfc8}\\
    &&&x^{ve}_{0} = 0;&& ~\forall v\in \mathcal{V}, e\in \mathcal{E}  \label{eq:mfc9}\\
    &&&\sum_{v\in \mathcal{V}} \sum_{e\in \mathcal{E}} x^{ve}_{i} \leq d_i;&&~\forall i\in N \label{eq:mfc10}\\
    &&&x^{ve}_{i} \leq q^{ve}_{i};&& ~\forall i\in N,v\in \mathcal{V},e\in \mathcal{E} \label{eq:mfc11}\\
    &&&\sum_{i\in N} x^{ve}_{i} \leq Q;&& ~\forall v\in \mathcal{V}, e\in \mathcal{E} \label{eq:mfc12}\\
    &&&q^{ve}_{j} \leq Q + M (1 - y^{ve}_{0j});&& ~\forall j\in N, v\in \mathcal{V}, e\in \mathcal{E} \label{eq:mfc13}\\
    &&&q^{ve}_{j} \geq Q - M (1 - y^{ve}_{0j});&& ~\forall j\in N, v\in \mathcal{V}, e\in \mathcal{E} \label{eq:mfc14}\\
    &&&q^{ve}_{j} \leq (q^{ve}_{i} - x^{ve}_{i}) + M (1 - y^{ve}_{ij});&& ~\forall i,j\in N,v\in \mathcal{V},e\in \mathcal{E} \label{eq:mfc15}\\
    &&&q^{ve}_{j} \geq (q^{ve}_{i} - x^{ve}_{i}) - M (1 - y^{ve}_{ij});&& ~\forall i,j\in N,v\in \mathcal{V},e\in \mathcal{E} \label{eq:mfc16}\\
    &&&q^{ve}_{i} <= M \lambda^{ve}_{i};&& ~\forall i\in N, v\in \mathcal{V}, e\in \mathcal{E} \label{eq:mfc17}\\
    &&&\underline{t}^{v0} = 0;&& ~\forall v\in \mathcal{V} \label{eq:mfc18}\\
    &&&\underline{t}^{v, e+1} = \underline{t}^{ve} + \sum_{i\in N}\sum_{j\in N} \tau_{ij}y^{ve}_{ij};&& ~\forall v\in \mathcal{V}, e\in \mathcal{E}\setminus 0  \label{eq:mfc19}\\
    &&&t^{ve}_{j} \geq \underline{t}^{ve} + \tau_{0j} - M(1-y^{ve}_{0j});&& ~\forall j\in N, v\in \mathcal{V}, e\in \mathcal{E} \label{eq:mfc20}\\
    &&&t^{ve}_{j} \geq t^{ve}_{i} + \tau_{ij} - M(1-y^{ve}_{ij});&& ~\forall i,j\in N\setminus0, v\in \mathcal{V}, e\in \mathcal{E}  \label{eq:mfc21}\\
    &&&t^{ve}_{j} \leq M\lambda^{ve}_{j};&& ~\forall j\in N, v\in \mathcal{V}, e\in \mathcal{E} \label{eq:mfc22}\\
    &&&\underline{t}^{v, E-1} + \sum_{i\in N}\sum_{j\in N} \tau_{ij}y^{v, E-1}_{ij} \leq L;&& ~\forall v\in \mathcal{V}  \label{eq:mfc23} \\
    &&& 0\leq x^{ve}_{i}, q^{ve}_{i} \leq Q, ~0\leq t^{ve}_{i},~ \underline{t}^{ve} \leq L, y^{ve}_{ij}, \lambda^{ve}_{i} \in \{0, 1\}; \qquad && \forall i,j\in N, v\in \mathcal{V}, e\in \mathcal{E} \nonumber
\end{align}
In this mathematical formulation, the objective function~\eqref{eq:mfc0} maximizes the total amount of served demands. Constraints \eqref{eq:mfc1}--\eqref{eq:mfc3} ensure that each customer cannot be visited more than once on each trip. Constraint \eqref{eq:mfc4} implies that if a vehicle visits no customer during a trip, no traverse between nodes should be planned for that vehicle on that trip. On the other hand, constraint \eqref{eq:mfc5} states that if a customer is decided to be served, a route to visit that customer should be planned. Constraints \eqref{eq:mfc6} and \eqref{eq:mfc7} guarantee that trips of a vehicle should be planned consecutively. Constraints \eqref{eq:mfc8}--\eqref{eq:mfc12} restrict the amount a vehicle can serve when visiting a customer on a trip by customer's demand and vehicle capacity. Constraints \eqref{eq:mfc13}--\eqref{eq:mfc17} are flow balance constraints. Lastly, constraints \eqref{eq:mfc18}--\eqref{eq:mfc23} keep track of operation time and restrict it by the given duration limit.
\begin{tiny}
\begin{algorithm}
\footnotesize
\caption{Compute an upper bound for the number of trips of a vehicle $E$}\label{alg:cap}
$C$ is the set of customers in $\mathcal{C}_0$ ordered such that $\tau_{0i} \leq \tau_{0j}, ~\forall i<j$\;
$E \gets 0$,
$i \gets 0$,
$l \gets 0$\;
\While{$l \leq L$}{
    \If{$l + \tau_{0i} + \tau_{i0} > L$}{
    \textbf{return} $E$\;
    }
    \Else{
        $l \gets l + \tau_{0i} + \tau_{i0}$\;
        $E \gets E+1$\;
        \If{$d_i \geq Q$}{
            $d_i \gets d_i - Q$\;
            }
        \Else{
            $i \gets i + 1$\;
            }
        }
}
\end{algorithm}
\end{tiny}

\section{Extended computational results}\label{app:exres}
This section provides extended results of the first and third experiment. Note that row Avg reports averages over instances solved to optimality by DAPP.
\begin{table}[!hbtp]
\footnotesize
\begin{tabular}{ccccccccccccc}
\hline
$\Psi_{n_z}$ & $Q$ & $n$ & $\mathbb{E}\sum d_c$ & DecQN & \%$\mathbb{E}\sum d_c$  & Opt. gap & DAPP & RP & GP & \%DAPP & \%RP & \%GP    \\ \hline
\multirow{22}{*}{\rotatebox[origin=c]{90}{Very Low}} & \multirow{11}{*}{25} & 12   & 125.0   & 74.6  & 59.7\% & 66.7\%   & -     & 53.3 & 55.5  & -      & 40.0\% & 34.3\% \\
                         &                      & 8    & 85.0    & 81.1  & 95.4\% & 0.0\%    & 75.4  & 58.3 & 56.6  & 7.6\%  & 39.0\% & 43.2\% \\
                         &                      & 9    & 100.0   & 61.8  & 61.8\% & 0.0\%    & 61.2  & 46.4 & 49.6  & 1.0\%  & 33.0\% & 24.5\% \\
                         &                      & 11   & 115.0   & 75.0  & 65.2\% & 0.0\%    & 76.6  & 51.2 & 52.2  & -2.1\% & 46.6\% & 43.5\% \\
                         &                      & 10   & 90.0    & 65.5  & 72.7\% & 0.0\%    & 61.0  & 49.1 & 56.0  & 7.3\%  & 33.4\% & 16.9\% \\
                         &                      & 11   & 95.0    & 52.8  & 55.5\% & 0.0\%    & 53.0  & 39.5 & 44.8  & -0.4\% & 33.7\% & 17.8\% \\
                         &                      & 12   & 125.0   & 81.7  & 65.4\% & 38.9\%   & -     & 59.7 & 69.2  & -      & 37.0\% & 18.1\% \\
                         &                      & 10   & 95.0    & 58.1  & 61.2\% & 0.0\%    & 55.8  & 44.5 & 51.4  & 4.1\%  & 30.6\% & 13.0\% \\
                         &                      & 9    & 115.0   & 74.0  & 64.4\% & 0.0\%    & 73.5  & 61.4 & 60.2  & 0.7\%  & 20.5\% & 22.9\% \\
                         &                      & 12   & 105.0   & 64.9  & 61.8\% & 24.0\%   & -     & 48.3 & 52.7  & -      & 34.3\% & 23.1\% \\ \cline{3-13}
                         &  \multicolumn{9}{r}{\textbf{Avg}} & \textbf{2.6\%}  & \textbf{33.8\%} & \textbf{26\%} \\ \cline{2-13}
                         & \multirow{11}{*}{50} & 12   & 125.0   & 102.4 & 81.9\% & 0.0\%    & 101.3 & 66.1 & 96.1  & 1.1\%  & 55.0\% & 6.5\%  \\
                         &                      & 8    & 85.0    & 84.8  & 99.8\% & 0.0\%    & 77.6  & 69.7 & 81.5  & 9.3\%  & 21.8\% & 4.1\%  \\
                         &                      & 9    & 100.0   & 85.3  & 85.3\% & 0.0\%    & 85.0  & 58.0 & 66.2  & 0.3\%  & 47.1\% & 28.8\% \\
                         &                      & 11   & 115.0   & 89.9  & 78.2\% & 0.0\%    & 89.0  & 61.1 & 64.1  & 1.1\%  & 47.1\% & 40.2\% \\
                         &                      & 10   & 90.0    & 85.8  & 95.3\% & 0.0\%    & 83.8  & 55.6 & 83.8  & 2.3\%  & 54.4\% & 2.3\%  \\
                         &                      & 11   & 95.0    & 65.3  & 68.8\% & 0.0\%    & 65.3  & 46.6 & 59.7  & 0.0\%  & 40.3\% & 9.4\%  \\
                         &                      & 12   & 125.0   & 108.8 & 87.0\% & 0.0\%    & 108.8 & 72.3 & 69.6  & 0.0\%  & 50.4\% & 56.3\% \\
                         &                      & 10   & 95.0    & 74.5  & 78.5\% & 0.0\%    & 74.5  & 50.5 & 74.5  & 0.0\%  & 47.6\% & 0.0\%  \\
                         &                      & 9    & 115.0   & 90.3  & 78.5\% & 0.0\%    & 89.7  & 71.9 & 80.4  & 0.7\%  & 25.6\% & 12.4\% \\
                         &                      & 12   & 105.0   & 85.8  & 81.7\% & 0.0\%    & 88.5  & 56.0 & 88.5  & -3.0\% & 53.3\% & -3.0\% \\ \cline{3-13}
                         & \multicolumn{9}{r}{\textbf{Avg}} & \textbf{1.2\%}  & \textbf{44.3\%} & \textbf{15.7\%} \\ \hline
\multirow{22}{*}{\rotatebox[origin=c]{90}{Low}}  & \multirow{11}{*}{25} & 14   & 145.0   & 105.7 & 72.9\% & 26.1\%   & -     & 67.5 & 87.7  & -      & 56.5\% & 20.5\% \\
                         &                      & 14   & 125.0   & 81.5  & 65.2\% & 38.9\%   & -     & 53.4 & 72.5  & -      & 52.7\% & 12.4\% \\
                         &                      & 18   & 165.0   & 112.9 & 68.4\% & 32.0\%   & -     & 71.0 & 101.1 & -      & 59.0\% & 11.6\% \\
                         &                      & 15   & 150.0   & 98.4  & 65.6\% & 42.9\%   & -     & 61.0 & 87.7  & -      & 61.2\% & 12.2\% \\
                         &                      & 10   & 100.0   & 82.7  & 82.7\% & 0.0\%    & 80.8  & 52.7 & 66.4  & 2.3\%  & 56.9\% & 24.5\% \\
                         &                      & 17   & 180.0   & 92.8  & 51.5\% & 50.0\%   & -     & 53.4 & 90.4  & -      & 73.6\% & 2.7\%  \\
                         &                      & 18   & 180.0   & 112.1 & 62.3\% & 50.0\%   & -     & 71.1 & 94.7  & -      & 57.6\% & 18.4\% \\
                         &                      & 18   & 175.0   & 101.3 & 57.9\% & 59.1\%   & -     & 61.6 & 87.9  & -      & 64.5\% & 15.3\% \\
                         &                      & 14   & 115.0   & 95.7  & 83.2\% & 9.5\%    & -     & 60.0 & 86.0  & -      & 59.6\% & 11.3\% \\
                         &                      & 12   & 140.0   & 113.2 & 80.9\% & 12.0\%   & -     & 76.5 & 99.0  & -      & 48.0\% & 14.3\% \\ \cline{3-13}
                         & \multicolumn{9}{r}{\textbf{Avg}} & \textbf{2.3\%}  & \textbf{56.9\%} & \textbf{24.5\%} \\ \cline{2-13}
                         & \multirow{11}{*}{50} & 14   & 145.0   & 132.3 & 91.2\% & 7.4\%    & -     & 82.6 & 103.6 & -      & 60.2\% & 27.7\% \\
                         &                      & 14   & 125.0   & 103.6 & 82.9\% & 19.1\%   & -     & 65.9 & 99.5  & -      & 57.2\% & 4.1\%  \\
                         &                      & 18   & 165.0   & 139.8 & 84.7\% & 10. 0\%  & -     & 88.1 & 111.6 & -      & 58.6\% & 25.3\% \\
                         &                      & 15   & 150.0   & 122.9 & 81.9\% & 20.0\%   & -     & 74.8 & 101.2 & -      & 64.3\% & 21.4\% \\
                         &                      & 10   & 100.0   & 96.2  & 96.2\% & 0.0\%    & 84.4  & 69.2 & 85.9  & 13.9\% & 39.1\% & 12.0\% \\
                         &                      & 17   & 180.0   & 117.6 & 65.3\% & 44.0\%   & -     & 78.1 & 105.1 & -      & 50.5\% & 11.9\% \\
                         &                      & 18   & 180.0   & 141.5 & 78.6\% & 24.1\%   & -     & 86.1 & 120.3 & -      & 64.4\% & 17.6\% \\
                         &                      & 18   & 175.0   & 119.4 & 68.3\% & 25.0\%   & -     & 77.7 & 112.6 & -      & 53.8\% & 6.1\%  \\
                         &                      & 14   & 115.0   & 110.2 & 95.8\% & 0.0\%    & 111.0 & 73.0 & 97.0  & -0.7\% & 51.1\% & 13.6\% \\
                         &                      & 12   & 140.0   & 136.2 & 97.3\% & 0.0\%    & 115.1 & 93.7 & 105.4 & 18.3\% & 45.4\% & 29.2\% \\ \cline{3-13}
                         & \multicolumn{9}{r}{\textbf{Avg}} & \textbf{10.5\%} & \textbf{45.2\%} & \textbf{18.2\%} \\ \hline
\multicolumn{10}{r}{\textbf{Total Avg}} & \textbf{3.0\%} & \textbf{41.5\%} & \textbf{19.9\%} \\ \hline
\end{tabular}
\caption{Extended results of DecQN compared with DAPP, RP, and GP on small instances}
\label{tab:resultsvrpscd0ex}
\end{table}

\begin{table}[!hbtp]
\centering
\footnotesize
\begin{tabular}{ccccccccc}
\hline
$L$ & $U$ & $\textrm{DecQN}_{L}$ & DecQN & GLQ & GHQ & \%DecQN & \%GLQ & \%GHQ \\ \hline
\multirow{3}{*}{S} & L & 765.6  & 799.7  & 738.3  & 834.0  & -4.25\% & 3.71\% & -8.19\% \\
                   & M & 737.1  & 781.0  & 718.5  & 794.7  & -5.62\% & 2.60\% & -7.24\% \\
                   & H & 722.4  & 751.0  & 704.9  & 760.4  & -3.81\% & 2.48\% & -5.00\% \\ \hline
\multirow{3}{*}{M} & L & 1027.0 & 1054.8 & 1019.2 & 1077.5 & -2.63\% & 0.77\% & -4.69\% \\
                   & M & 1017.7 & 1035.8 & 1002.8 & 1067.0 & -1.75\% & 1.49\% & -4.62\% \\
                   & H & 1002.8 & 1023.8 & 977.1  & 1043.6 & -2.06\% & 2.63\% & -3.91\% \\ \hline
\multirow{3}{*}{L} & L & 1079.0 & 1078.7 & 1078.3 & 1078.3 & 0.02\%  & 0.06\% & 0.06\%  \\
                   & M & 1078.2 & 1080.3 & 1076.3 & 1076.9 & -0.19\% & 0.18\% & 0.12\%  \\
                   & H & 1079.8 & 1081.4 & 1072.3 & 1075.3 & -0.15\% & 0.70\% & 0.42\%  \\ \hline
                   &   &        &        &        & \textbf{Avg}   & \textbf{-2.3\%}  & \textbf{1.6\%}  & \textbf{-3.7\%}  \\ \hline
\end{tabular}
\caption{Extended results of $\textrm{DecQN}_L$}
\label{tab:exaggrdl}
\end{table}

\begin{table}[!hbtp]
\centering
\footnotesize
\begin{tabular}{ccccccccc}
\hline
$L$ & $U$ & $\textrm{DecQN}_{U}$ & DecQN  & GLQ & GHQ & \%DecQN & \%GLQ & \%GHQ \\ \hline
\multirow{3}{*}{S} & L & 773.5  & 799.7  & 738.3  & 834.0  & -3.3\%  & 4.8\% & -7.2\% \\
                   & M & 760.2  & 781.0  & 718.5  & 794.7  & -2.7\%  & 5.8\% & -4.3\% \\
                   & H & 745.3  & 751.0  & 704.9  & 760.4  & -0.8\%  & 5.7\% & -2.0\% \\ \hline
\multirow{3}{*}{M} & L & 1051.7 & 1054.8 & 1019.2 & 1077.5 & -0.3\%  & 3.2\% & -2.4\% \\
                   & M & 1038.6 & 1035.8 & 1002.8 & 1067.0 & 0.3\%   & 3.6\% & -2.7\% \\
                   & H & 1024.6 & 1023.8 & 977.1  & 1043.6 & 0.1\%   & 4.9\% & -1.8\% \\ \hline
\multirow{3}{*}{L} & L & 1078.7 & 1078.7 & 1078.3 & 1078.3 & 0.0\%   & 0.0\% & 0.0\%  \\
                   & M & 1080.3 & 1080.3 & 1076.3 & 1076.9 & 0.0\%   & 0.4\% & 0.3\%  \\
                   & H & 1080.3 & 1081.4 & 1072.3 & 1075.3 & -0.1\%  & 0.7\% & 0.5\%  \\ \hline
                   &   &        &        &        & \textbf{Avg}   & \textbf{-0.7\%}  & \textbf{3.2\%} & \textbf{-2.2\%} \\ \hline
\end{tabular}
\caption{Extended results of $\textrm{DecQN}_U$}
\label{tab:exresultsaggrsv}
\end{table}

\begin{table}[!hbtp]
\centering
\footnotesize
\begin{tabular}{ccccccccc}
\hline
$L$ & $U$ & $\textrm{DecQN}_{All}$ & DecQN  & GLQ & GHQ & \%DecQN & \%GLQ & \%GHQ \\ \hline
\multirow{3}{*}{S} & L & 755.9  & 799.7  & 738.3  & 834.0  & -5.5\% & 2.4\%  & -9.4\% \\
                   & M & 743.0  & 781.0  & 718.5  & 794.7  & -4.9\% & 3.4\%  & -6.5\% \\
                   & H & 731.6  & 751.0  & 704.9  & 760.4  & -2.6\% & 3.8\%  & -3.8\% \\ \hline
\multirow{3}{*}{M} & L & 1026.5 & 1054.8 & 1019.2 & 1077.5 & -2.7\% & 0.7\%  & -4.7\% \\
                   & M & 1019.2 & 1035.8 & 1002.8 & 1067.0 & -1.6\% & 1.6\%  & -4.5\% \\
                   & H & 1009.3 & 1023.8 & 977.1  & 1043.6 & -1.4\% & 3.3\%  & -3.3\% \\ \hline
\multirow{3}{*}{L} & L & 1078.1 & 1078.7 & 1078.3 & 1078.3 & -0.1\% & 0.0\% & 0.0\% \\
                   & M & 1078.7 & 1080.3 & 1076.3 & 1076.9 & -0.1\% & 0.2\%  & 0.2\%  \\
                   & H & 1077.9 & 1081.4 & 1072.3 & 1075.3 & -0.3\% & 0.5\%  & 0.2\%  \\ \hline
                   &   &        &        &        & \textbf{Avg }   & \textbf{-2.1\%} & \textbf{1.8\%}  & \textbf{-3.5\%} \\ \hline
\end{tabular}
\caption{Extended results of $\textrm{DecQN}_{All}$}
\label{tab:exresultsaggrall}
\end{table}

\section{On-line-off-line DecQN}\label{sec:resultsonoff}
In the previous experiments, we investigated the performance of the trained Q-network as a policy. More specifically, for a given Q-network that provides Q values for observation-action pairs, the policy applied the action with the maximum Q value. This section investigates the use of the trained Q-network as a reward-to-go approximation function in an on-line RA. The approach presented in this experiment is similar to the solution method proposed in~\cite{Ulmer2017}, where an off-line trained value function is embedded into an RA in order to approximate the reward-to-go.
In Section \ref{sec:resultsonoffra}, the proposed RA is explained in detail. We use instances and trained policies as illustrated in Section 6.1. Finally, we compare the performance of the proposed on-line-off-line DecQN with the purely off-line method in Section \ref{sec:resultsonoffresults}.

\subsection{Proposed Rollout Algorithm}\label{sec:resultsonoffra}
The Rollout Algorithm is an on-line method in the form of forward-dynamic programming that can be viewed as a one-step policy iteration \citep{Bertsekas2013}. 
\cite{Powell2011} showed that if we roll out a known heuristic policy (the so-called base policy) to approximate reward-to-go values, the resulting policy will be improved over the base policy. In practice, the rollout algorithm estimates the reward-to-go value of an available action at the current decision epoch by simulating the base policy after choosing that action. Interested readers may refer to a survey by \cite{Bertsekas2013} for a detailed discussion about RAs.

Similar to what defined in Equation (20), we let $V^{\pi}(s_k)$ be the expected reward that can be obtained by starting from the state $s_{k}$ and following the base policy $\pi$. Accordingly, we define the rollout policy $\Pi$ as:
\begin{equation} \label{eq:rollout}
    \Pi(s_{k}) = \arg\max_{x_k\in A(s_k)} \left[R(s_k, x_k) +  \mathop{\mathbb{E}}_{w\in W}[\gamma V^{\pi}(s_{k+1})]\right].
\end{equation}
In our case, the base policy is provided by the Q values trained by the DecQN method in the previous experiment.
Therefore, we substitute the approximate value function $V^{\pi}(s_{k+1})$ according to Equation (32). Also, following the partially decentralized formulation, we replace $s_k$ and $A(s_k)$ with $o_{k, \bar{v}}$ and $A(s_k, \bar{v})$, respectively. We rewrite the rollout policy $\Pi$ as:
\begin{equation} \label{eq:rollout2}
    \Pi(o_{k,\bar{v}}) = \arg\max_{x_k\in A(s_k, \bar{v})} \left[R(o_{k, \bar{v}}, x_k) +  \mathop{\mathbb{E}}_{w\in W}[\gamma \max_{x'\in A(s_{k+1}, \bar{v}')} Q(o_{k+1, \bar{v}'}, x')]\right],
\end{equation}
where $\bar{v}'$ represents the active vehicle in decision epoch $k+1$.

Note that, since the set of customers is realized at the beginning of the operation, the customer demands are the only source of uncertainty in the equation above. Therefore, the expectation function $\mathbb{E}(.)$ is computed over a finite set of sample scenarios for the realization of the stochastic demands ($W$). It is important to note that, according to Equation (18), which defines the reward function for DecQN, the actual reward is realized with a delay. However, it is not feasible to wait to receive the actual reward through on-line methods, such as RA. 
Therefore, instead of using the reward function as in Equation (18), we define it as follows:
\begin{equation}
    R(o_{k, \bar{v}}, x_k) = \begin{cases}
    \Tilde{d}_c & x_k \textrm{ is a } c\in \mathcal{C}_0 \\
    0 & \textrm{otherwise}
    \end{cases},
\end{equation}
where $\Tilde{d}_c$ is the unserved demand of customer $c$, if its actual demand is already realized; otherwise, it is the expected demand.

The RA proceeds at decision epoch $k$ by simulating the system for each possible action $x_k\in A(s_k, \Bar{v})$ for one step and observe the immediate reward $R(o_{k, \bar{v}}, x_k)$ as well as the next observation $o_{k+1,\bar{v}'}$. In the next observation, the maximum Q value of available actions is considered as the approximate reward expected to be obtained onward, denoted as $\bar{f}$. The expected reward-to-go of each action is computed for a set of demand realizations ($W$). The action with the maximum $\Bar{f}$ will be selected as the best action to take ($x^*$) at decision epoch $k$.
\begin{algorithm}[!hbtp]
\footnotesize
 \SetAlgoLined
 \DontPrintSemicolon
 \LinesNumbered
 Initialize: set $f^*$ to 0 \;
 Construct the action space $A(s_k, \Bar{v})$\;
 \For{$x_k\in A(s_k, \Bar{v})$}
 {
    $\Bar{f} \gets 0$ \;
    $o_{k, \Bar{v}} \gets O(s_k, \Bar{v})$\;
    Take action $x_k$ and receive the reward $r_k=R(o_{k, \bar{v}}, x_k)$\;
    \For{$w \in W$}{
        State transition $k$ to $k+1$:\;
        $s_{k+1} \gets S^M(s_k, x_k, w)$ \;
        $\Bar{v}' \gets \{v\in\mathcal{V}|a_v=\min_{v \in \mathcal{V}} a_v\}$ \;
        $o_{k+1, \Bar{v}'} \gets O(s_{k+1}, \Bar{v}') $\;
        Take action $x_{k+1} = \arg\max_{x'} Q(o_{k+1, \Bar{v}'}, x')$ and receive $r'_k=R(o_{k+1, \bar{v}'}, x_{k+1})$\;
        state transition $k+1$ to $k+2$:\;
        $s_{k+2} \gets S^M(s_{k+1}, x_{k+1}, w)$ \;
        $\Bar{v}'' \gets \{v\in\mathcal{V}|a_v=\min_{v \in \mathcal{V}} a_v\}$ \;
        $o_{k+2, \Bar{v}''} \gets O(s_{k+2}, \Bar{v}'') $\;
        $f \gets [r_k + \gamma r'_k + \gamma^2\max_{x_{k+2}\in A(s_{k+2}, \Bar{v}'')} Q(o_{k+2, \Bar{v}''}, x_{k+2})]$\;
        $\Bar{f} \gets \Bar{f} + f $ \;
    }
    $\Bar{f} \gets \Bar{f} / |W|$\;
    \If{$\Bar{f} > f^*$}
    {
        $x^* \gets x_k$\;
        $f^* \gets \Bar{f}$\;
    }
 }
 \textbf{return} $x^*$\;
 \caption{Proposed Rollout Algorithm for decision epoch $k$}
 \label{alg:ra}
\end{algorithm}

Given the delayed reward mechanism, we decided to simulate the problem for two steps for each action (instead of the traditional one-step) in our implementation. Therefore, the implemented Equation \eqref{eq:rollout2} has been modified as:
\begin{equation} \label{eq:rollout3}
    \Pi(o_{k,\Bar{v}}) = \arg\max_{x_k\in A(s_k,\Bar{v})} \mathop{{}\mathbb{E}}_{w\in W}[R(o_{k, \bar{v}}, x_k) + \gamma R(o_{k+1, \bar{v}'}, x_{k+1}) + \gamma^2 \max_{x_{k+2}\in A(s_{k+2},\Bar{v}'')} Q(o_{k+2,\Bar{v}''}, x_{k+2})],
\end{equation}
where $x_{k+1}$ is the action in decision epoch $k+1$ according to the base policy $\pi$ (i.e., $\displaystyle{x_{k+1}= \arg\max_{x'\in A(s_k+1,\Bar{v}')} Q(o_{k+1,\Bar{v}'}, x')}$) and $\Bar{v}''$ is the active vehicle in decision epoch $k+2$. Algorithm \ref{alg:ra} illustrates our proposed RA for each decision epoch $k$, when the active vehicle is $\bar{v}$, the global state is $s_k$, and $W$ is a finite set of demand realizations.

\subsection{Results and Discussion} \label{sec:resultsonoffresults}
In Section 6.1, we tested each trained policy for 250,000 realizations. However, because the rollout algorithm requires substantial on-line computing, testing across the same realization set would be computationally expensive. Therefore, we test the proposed on-line-off-line DecQN for each instance $i\in I$, with $I=\mathcal{D}\times \mathcal{Q}$, on $\hat{I}(i)$ with 2,500 realizations (50 customer realizations $\times$ 50 demand realizations). Furthermore, we only focus on harder instances with $\mathcal{D}=\{\textrm{Moderate, High, Very High}\}$. For this experiment, we set the size of the sample scenarios for the rollout algorithm ($|W|$) to 50.
\begin{table}[!hbtp]
\centering
\footnotesize
\begin{tabular}{ccccc}
\hline
$\Psi_{n_z}$ & $Q$ & RA-DecQN & DecQN& \%DecQN \\ \hline
\multirow{3}{*}{Moderate}  & 25 & 164.3  & 161.3 & 1.9\%  \\
                    & 50 & 202.7 &	196.9  & 3.0\% \\
                    & 75 & 217.1 &	211.9 & 2.5\%  \\ \hline
\multirow{3}{*}{High}  & 25 & 358.3 &	358.2 &  0.0\%  \\
                    & 50 & 461.0 &	455.4 & 1.2\%  \\
                    & 75 & 510.5 &	501.5 & 1.8\%  \\ \hline
\multirow{3}{*}{Very High}  & 25 & 551.7 &	551.7 &	0.0\% \\
                    & 50 & 708.7 &	701.4 &	1.0\%  \\
                    & 75 & 790.3 & 776.7 & 1.8\%  \\ \hline
\multicolumn{4}{r}{\textbf{Avg}} & \textbf{1.5\%} \\ \hline
\end{tabular}
\caption{Results of on-line-off-line DecQN compared with off-line DecQN}
\label{tab:raresults}
\end{table}
Table \ref{tab:raresults} reports the results of the on-line-off-line DecQN (RA-DecQN) and compares them with the performance of the off-line policy (DecQN). The average improvement of the rollout algorithm over the off-line DecQN is seen in column \%DecQN. According to the results, the proposed rollout outperforms the off-line version by an average of 1.5\% in all cases. 
The average improvement of the rollout algorithm over the off-line DecQN is shown in Table \ref{tab:resultsra2}. This table helps us to analyze how changes in customer density and vehicle capacity affect \%DecQN. Accordingly, as the customer density decreases from Very High to Moderate, the gap widens significantly. Specifically, the average improvement over the off-line DecQN in instances with Very High customer density is 0.9\%, whereas this value in instances with Moderate density is 2.4\%. It is essential to notice that, as discussed in Section 6.1, the importance of a good routing strategy is magnified in instances with a lower customer density. This finding shows that the on-line-off-line method improves the off-line policy in terms of routing aspects.
Furthermore, the sensitivity analysis over the changes in the vehicle capacity also confirms this finding. According to Table \ref{tab:resultsra2}, the improvement percentage considerably increases from 0.6\% to 2.0\% when the vehicle capacity changes from 25 to 75. It is worth to notice that as the vehicle capacity becomes larger, vehicles are allowed to perform longer routes which highlights the necessity of an effective routing strategy. 
\begin{table}[!hbtp]
\centering
\footnotesize
\begin{tabular}{cccccc}
\hline
\multicolumn{2}{c}{\multirow{2}{*}{\%DecQN}} & \multicolumn{3}{c}{$Q$} & \\ \cline{3-5}
\multicolumn{2}{c}{}  & 25 & 50     & 75     & \textbf{Avg}    \\ \hline
\multirow{3}{*}{$\Psi_{n_z}$} & Moderate   & 1.9\% & 3.0\% & 2.5\% & \textbf{2.4\%} \\
                   & High   & 0.0\% & 1.2\% & 1.8\% & \textbf{1.0\%} \\
                   & Very High   & 0.0\% & 1.0\% & 1.8\% & \textbf{0.9\% }\\ \hline
 \multicolumn{2}{r}{\textbf{Avg}} & \textbf{0.6\%} & \textbf{1.7\%} & \textbf{2.0\%} & \textbf{1.5\%} \\  \hline
\end{tabular}
\caption{Average improvement of on-line-off-line DecQN over off-line DecQN}
\label{tab:resultsra2}
\end{table}
In Section 6.1, we demonstrated that the average improvement over GP drops from 11.2\% to 7.9\%, when the vehicle capacity increases from 25 to 75. However, this experiment shows that the on-line-off-line DecQN partly compensates for this reduction.

In this experiment, the average computing time per decision epoch for Moderate, High, and Very High customer density is 0.10, 0.17, and 0.21 seconds, respectively. Although, as expected, the computation time grows as the number of customers increases, the on-line-off-line DecQN can decide in 0.21 seconds, even in the most complex problems. Therefore, the proposed rollout algorithm can also be regarded as an efficient on-line method.